\documentclass[11pt]{article}

\usepackage{amsmath,amsthm,amssymb}
\usepackage{mathtools}
\usepackage[bb=boondox]{mathalfa}
\usepackage{dsfont}
\usepackage{booktabs}
\usepackage{tabularx}
\usepackage{colonequals}
\usepackage{authblk}

\usepackage{graphicx}
\usepackage{tikz}
\usetikzlibrary{arrows.meta,positioning,shapes.geometric,backgrounds,fit,calc}
\usepackage{color}
\usepackage{soul}
\usepackage[mathlines]{lineno}
\usepackage[numbers]{natbib}
\usepackage{centernot}
\usepackage{float}
\usepackage{hyperref}
\usepackage[capitalise,nameinlink]{cleveref}

\usepackage[letterpaper,margin=1in]{geometry}
\renewcommand{\arraystretch}{1.2}

\numberwithin{figure}{section}
\numberwithin{table}{section}
\numberwithin{equation}{section}

\theoremstyle{plain}

\newtheorem{ctr}{Counter}[section]
\newtheorem{lem}[ctr]{Lemma}
\newtheorem{prop}[ctr]{Proposition}
\newtheorem{cor}[ctr]{Corollary}

\newtheorem{defn}[ctr]{Definition}
\newtheorem{rem}[ctr]{Remark}
\newtheorem{thm}[ctr]{Theorem}

\crefname{lem}{Lemma}{Lemmas}
\Crefname{lem}{Lemma}{Lemmas}
\crefname{prop}{Proposition}{Propositions}
\Crefname{prop}{Proposition}{Propositions}
\crefname{cor}{Corollary}{Corollaries}
\Crefname{cor}{Corollary}{Corollaries}
\crefname{defn}{Definition}{Definitions}
\Crefname{defn}{Definition}{Definitions}
\crefname{rem}{Remark}{Remarks}
\Crefname{rem}{Remark}{Remarks}
\crefname{thm}{Theorem}{Theorems}
\Crefname{thm}{Theorem}{Theorems}

  \ifdefined\AddToHook
  \AddToHook{env/lem/begin}{\crefalias{ctr}{lem}}
  \AddToHook{env/prop/begin}{\crefalias{ctr}{prop}}
  \AddToHook{env/cor/begin}{\crefalias{ctr}{cor}}
  \AddToHook{env/obs/begin}{\crefalias{ctr}{obs}}
  \AddToHook{env/fact/begin}{\crefalias{ctr}{fact}}
  \AddToHook{env/conj/begin}{\crefalias{ctr}{conj}}
  \AddToHook{env/quest/begin}{\crefalias{ctr}{quest}}
  \AddToHook{env/defn/begin}{\crefalias{ctr}{defn}}
  \AddToHook{env/rem/begin}{\crefalias{ctr}{rem}}
  \AddToHook{env/thm/begin}{\crefalias{ctr}{thm}}
  \fi

\definecolor{oiBlack}{HTML}{000000}
\definecolor{red}{rgb}{.8,0,0}     
\definecolor{bblu}{rgb}{0,0,1}     
\definecolor{gray}{rgb}{.4,.4,.4}  
\definecolor{green}{rgb}{0,.7,0}   

\newcommand{\Z}{{\mathbb Z}}
\newcommand{\R}{{\mathbb R}}

\newcommand{\bone}{\mathbb{1}}

\newcommand{\tw}{\operatorname{tw}}

\DeclareMathOperator{\gex}{gex}

\newcommand{\WL}[1]{\equiv_{\mathrm{WL}(#1)}}

\newcommand{\WLp}[1]{\mathrm{WL}(#1)}

\newcommand{\T}{\mathcal{T}}
\newcommand{\F}{\boldsymbol{\mathcal F}}
\newcommand{\LL}{\boldsymbol{\mathcal L}}
\newcommand{\G}{\mathcal G}

\newcommand{\Lr}{\Lambda_r}
\newcommand{\LrS}{\Lambda_{r,S}}

\newcommand{\noi}{\noindent}

\title{Transfer of Weisfeiler--Leman equivalence under graph operations}
\author{Dmitriy Kunisky\thanks{Email: \texttt{kunisky@jhu.edu}.}\,\,}
\author{Hamilton Sawczuk\thanks{Email: \texttt{hsawczu1@jhu.edu}.}}
\affil{Department of Applied Mathematics \& Statistics, Johns Hopkins University}
\date{\today}

\begin{document}
\pagenumbering{roman}
\pagestyle{empty}
\maketitle
\thispagestyle{empty}

\begin{abstract}
The Weisfeiler--Leman algorithms $\WLp{k}$ for $k\ge1$ form a hierarchy of increasingly
powerful tests for graph isomorphism. We study how the equivalence relations they induce on graphs
transfer through a graph operation $\Phi$: for a given $k$, what level of WL-equivalence
between graphs $H_1$ and $H_2$ guarantees $\WLp{k}$-equivalence of graphs $\Phi(H_1)$ and
$\Phi(H_2)$? We call the least such level the \emph{transfer parameter} $\tau_\Phi(k)$ and
determine it exactly for the line graph operation $\Phi = L$: $\tau_L(k)=2k+1$ for $k\ge2$, while
$\tau_L(1)=1$. We also extend these results to the $(r,S)$-line graphs, which generalize the line
graph from edges to $r$-cliques by connecting two cliques when their intersection size lies
in a prescribed set $S$. For these operations we prove that $\tau_{r,S}(1)=r-1$ for every $S$ and, for certain choices of $S$, that $\tau_{r,S}(k)=rk+r-1$ for $k\ge r$,
and we bound the transfer parameter in the remaining range. Our upper bounds on transfer parameters follow
from a treewidth analysis of diagrammatic expansions of homomorphism counts into $\Phi(H)$ as linear
combinations of homomorphism counts into $H$. The matching lower bounds combine explicit
constructions of witness graphs with the homomorphism-distinguishing closedness of graphs of bounded
treewidth.
\end{abstract}

\clearpage

\tableofcontents

\clearpage
\pagenumbering{arabic}
\pagestyle{plain}

\section{Introduction}
\label{sec:intro}

We write $\mathcal{G}$ for the set of all isomorphism classes of finite simple graphs with vertex set $[n]$ for some $n\geq0$.
Without further qualification all graphs we discuss below are taken to be simple, i.e., belonging to $\mathcal G$.
\emph{Graph operations} are functions $\Phi: \mathcal{G} \to \mathcal{G}$ that take an input graph and generate another graph from it.
One of the most widely-studied graph operations is the formation of the \emph{line graph}:

\begin{defn}[Line graph]
\label{defn:lg}
    The \emph{line graph} $L(H)$ of a graph $H = (V, E)$ has the edge set $E$ as its vertex set,
    with pairs of edges adjacent if they share an endpoint. Formally,
    \[
        V(L(H))=E,\qquad E(L(H))=\{\{e_1, e_2\}  : e_1,e_2\in E,\;|e_1\cap e_2|=1\}.
    \]
\end{defn}
\noindent
Many graph properties transform conveniently under the line graph operation; for instance, a
matching in $H$ is an independent set in $L(H)$, and an edge-coloring in $H$ is a
vertex-coloring in $L(H)$.
Classically, these correspondences provide reductions among combinatorial problems on graphs: by running a vertex-coloring algorithm on $L(H)$, we can find edge-colorings of $H$, and likewise for other problems~\cite{West2001}.

Graph operations are used, for example, in graph learning, where the output of graph operations such as the line graph can be used as ``features.''
To learn a function on graphs $H \mapsto f(H) \in \mathbb{R}$, one can instead learn a function of $L(H)$ or of $H$ and $L(H)$ jointly.
In the simplest case, the power of an algorithm can be ``boosted'' by running it on a graph generated by a graph operation.
For example, message-passing neural networks are limited by the expressive power of the
Weisfeiler--Leman algorithm discussed below~\cite{Morris19,Xu19}, but running an analogous network on $L(H)$ rather than
$H$ can distinguish inputs the original cannot, as in the line-graph neural network of Chen,
Li, and Bruna~\cite{ChenLiBruna19} and the edge-reified message passing of Liu and
Song~\cite{LiuSong22}.

We consider the interaction of the line graph with a common technique for refuting \emph{graph isomorphism}, the \emph{Weisfeiler--Leman (WL) hierarchy} of algorithms.
Given two input graphs $H_1, H_2$, these algorithms can be viewed as outputting either ``not isomorphic'' or ``no result'', where the former output means the algorithm has found a certificate that $H_1 \not\cong H_2$, while the latter output means it could be that either $H_1 \cong H_2$ or $H_1 \not\cong H_2$.
These algorithms come in a hierarchy of increasing power and computational complexity, indexed by a parameter $k \geq 1$, the $k$th of which we denote $\WLp{k}$; in particular  if $\WLp{k}$ outputs ``not isomorphic,'' then so does $\WLp{k^{\prime}}$ for each $k^{\prime} > k$.

The WL algorithms can be instantiated via iterative coloring~\cite{CFI92,GroheNeuen21}, convex optimization~\cite{GRS22,RobersonSeppelt23}, and graph neural
networks~\cite{Morris19,Xu19}. They are also characterized more abstractly by homomorphism
counts from bounded-treewidth graphs~\cite{Dvorak10,DGR18}, bounded-variable counting
logic~\cite{CFI92}, and pebbling comonads~\cite{ADW17,DJR21}.
We give the homomorphism count interpretation, which will be most suitable for our purposes.

\begin{defn}[Weisfeiler--Leman equivalence]
\label{defn:wl} For $k\ge1$ and graphs $H_1,H_2$, we say $H_1$ and $H_2$ are \emph{$\WLp{k}$-equivalent} and write $H_1\WL{k}H_2$ if
    \[
        \hom(G,H_1)=\hom(G,H_2)
    \]
    for every graph $G$ of treewidth at most $k$ (\cref{defn:tree-decomp}).
    Here $\hom(G, H)$ denotes the number of \emph{graph homomorphisms} of $G$ into $H$, i.e. functions $V(G)\to V(H)$ mapping edges to edges.
\end{defn}
\noindent
It follows that the relations $\equiv_{\WLp{k}}$ are equivalence relations on graphs, and these equivalence relations become finer as $k$ increases, so that whenever $k \leq k^{\prime}$, $H_1 \equiv_{\WLp{k^{\prime}}} H_2$ implies $H_1 \equiv_{\WLp{k}} H_2$.

We investigate how WL equivalence interacts with graph operations.
Namely, for a graph operation $\Phi: \mathcal{G} \to \mathcal{G}$ and a given $k \geq 1$, we ask what level of WL equivalence of $H_1$ and $H_2$ is needed to
guarantee $\WLp{k}$ equivalence of $\Phi(H_1)$ and $\Phi(H_2)$.
That is, what is the minimum value of $t$ such that
\[
    H_1\WL{t}H_2 \,\,\,\, \Longrightarrow \,\,\,\, \Phi(H_1)\WL{k}\Phi(H_2)\,\,\text{?}
\]
We refer to this question as the \emph{transfer of WL equivalence} through the graph operation $\Phi$.
Below we will answer this question exactly for $\Phi = L$, and give partial results for a family of generalizations.

The practical significance of the transfer parameter is algorithmic. While no known algorithm solves graph isomorphism in time polynomial in $n$~\cite{GroheNeuen21}, for any fixed $k$ the $\WLp{k}$ algorithm runs in $n^{O(k)}$ time, and in an average-case sense $\WLp{1}$ already distinguishes almost all pairs of graphs~\cite{BES80}.
Because of this, $\WLp{1}$ is a workhorse method for practical graph isomorphism, powering isomorphism solvers~\cite{McKayPiperno14} and
appearing as a subroutine in Babai's state-of-the-art quasipolynomial-time graph isomorphism
algorithm~\cite{Babai16}.
However, there are simple examples of non-isomorphic graphs that $\WLp{1}$ fails to distinguish, for example any pair of $d$-regular graphs on the same number of vertices.
These limitations can be addressed by passing to more expensive $\WLp{k}$ algorithms, but one may apply the above graph learning idea and instead run $\WLp{k}$ on the transformed graphs $\Phi(H_1)$ and $\Phi(H_2)$. The transfer parameter bounds the effectiveness of this approach in the worst case. More abstractly, WL equivalence can be viewed as a quantitative refinement of graph isomorphism.
For instance, we may define the ``WL similarity'' between graphs as $\Delta_{\mathrm{WL}}(H_1, H_2) \colonequals \max\{k: H_1 \equiv_{\WLp{k}} H_2\}$.
The transfer parameter then bounds how much this similarity can decrease under a given graph operation.

\subsection{Main results}

Our first main result identifies the transfer parameter for the line graph.
\begin{thm}[Transfer of WL-equivalence for the line graph]\label{thm:main-lg}
    Let $\tau_L(k)$ denote the minimum value of $t$ such that
    \[
        H_1\WL{t}H_2\;\Longrightarrow\;L(H_1)\WL{k}L(H_2)
    \]
    for all $H_1, H_2 \in \mathcal{G}$.
    Then, we have
    \begin{equation}\label{eqn:main-lg}
        \tau_L(k)=\begin{cases}
            1 & \text{if } k=1, \\
            2k+1 & \text{if } k\geq 2. \\
        \end{cases}
    \end{equation}
\end{thm}
\noindent
For connected graphs the case $k=1$ also follows from covering theory (\cref{rem:covering}); we include a proof using our methods in \cref{sec:k1}. For $k\geq2$, the upper bound $\tau_L(k)\le2k+1$ follows from finite-variable counting logic (\cref{subsec:related-work}), but the matching lower bound is new (\cref{sec:lower}). We partially extend \cref{thm:main-lg} to a variant of the line graph based on overlapping cliques, as follows.

\begin{defn}[$(r,S)$-line graph]\label{defn:rs}
    Fix an integer $r\ge2$ and a set $S\subseteq\{0,1,\dots,r-1\}$. For a graph $H$, the
    \emph{$(r,S)$-line graph} $\LrS(H)$ has the $r$-cliques of $H$ as its vertices, with
    two cliques $P,Q$ adjacent whenever $|P\cap Q|\in S$.
\end{defn}
\noindent
The vertices of $\LrS(H)$ are the $(r-1)$-simplices of the clique complex of $H$, and $S$ specifies those intersection sizes that induce adjacency.
\Cref{tab:rs-examples} lists three familiar special cases.
On $K_n$, the adjacency matrix of $\LrS(K_n)$ is $\sum_{s\in S}A_{r-s}$, where $A_j$ is the Johnson scheme matrix joining $r$-subsets with intersection size $r-j$~\cite{BannaiIto84,BCN89}.
Equivalently, this is the generalized Kneser graph with prescribed intersection set $S$~\cite{Haviv18}.

\begin{table}
\centering
\renewcommand{\arraystretch}{1.3}
\begin{tabularx}{\linewidth}{@{}l >{\raggedright\arraybackslash}X >{\raggedright\arraybackslash}X@{}}
    \toprule
    $(r,S)$ & Operation & Output on $K_n$ \\
    \midrule
    $(2,\{1\})$ & line graph & $J(n,2,1)$ (Johnson graph) \\
    $(r,\{r-1\})$ & $r$-line graph (triangle graph for $r=3$)~\cite{LakshmananBujtasTuza,GeneralizedLineGraphs15} & $J(n,r,r-1)$~\cite{ChenLih87} \\
    $(r,\{1,\hdots,r-1\})$ & line graph of $r$-uniform hypergraph of $r$-cliques~\cite{Berge73,NRSS82} & complement of the Kneser graph, $\overline{K(n,r)}$ \\
    \bottomrule
\end{tabularx}
\caption{Special cases of the $(r,S)$-line graph and their outputs on the complete graph $K_n$.
Here the generalized Johnson graph $J(n,r,s)$ joins $r$-subsets of $[n]$ whose intersection has size $s$.}
\label{tab:rs-examples}
\end{table}

\begin{thm}[Transfer of WL-equivalence for the $(r,S)$-line graph]\label{thm:main-rs}
    Fix $r\ge2$ and $S\subseteq\{0,1,\dots,r-1\}$. Let $\tau_{r,S}(k)$ denote the minimum value of $t$ such that
    \[
        H_1\WL{t}H_2\;\Longrightarrow\;\Lambda_{r,S}(H_1)\WL{k}\Lambda_{r,S}(H_2)
    \]
    for all $H_1,H_2 \in \mathcal{G}$.
    The following then hold:
    \begin{enumerate}
        \item $\tau_{r,S}(1) = r - 1$, and $r-1\leq\tau_{r,S}(k)\leq rk+r-1$ for every $k\geq1$.
        \item If $\{1\}\subseteq S\subseteq\{1,\dots,r-1\}$ and $k\geq r$, then $\tau_{r,S}(k) = rk + r - 1$.
        \item If $\{1\}\subseteq S\subseteq\{1,\dots,r-1\}$, $r\geq3$, and the sign condition \eqref{eqn:dagger} holds, then $\tau_{r,S}(k)\geq 2r-2$ for $2\leq k\leq r-1$.
    \end{enumerate}
\end{thm}
\noindent
Note that letting $r=2,S=\{1\}$ recovers \cref{thm:main-lg}. The sign condition \eqref{eqn:dagger} concerns the integers
\[
    \alpha_t=\sum_{\substack{s\in S\\s\le t}}(-1)^{t-s}\binom ts\qquad(0\le t\le r),
\]
which are the coefficients of the gadget expansion of $\Lambda_{r,S}$ (\cref{lem:coefficients}); it requires that $(-1)^{t-1}\alpha_t\geq0$ for $1\le t\le r$ and that $\alpha_{r-1}\neq0$. It holds, for example, when $S=\{1\}$ and when $S=\{1,\dots,r-1\}$ (\cref{lem:dagger-suff}).

We leave the problem of determining $\tau_{r,S}(k)$ for $2\le k\le r-1$ open; already for $r=3$ and $S=\{1\}$, our bounds give only $4\le\tau_{3,\{1\}}(2)\le8$.
Another area for future work is the case $S=\{r-1\}$ for $r\ge3$, which includes the triangle graph operation and does not satisfy the hypothesis $1\in S$ of items 2 and 3.

\subsection{Proof techniques: Diagrammatic expansion of graph operations}
\label{sec:techniques}

We first consider the task of proving upper bounds on the transfer parameter of a graph operation $\Phi$.
Unraveling the definitions, this requires us to show that whenever $\hom(G, H_1) = \hom(G, H_2)$ for all $G$ of treewidth at most some $t$, $\hom(G, \Phi(H_1)) = \hom(G, \Phi(H_2))$ for all $G$ of treewidth at most some $k$.
Our analysis proceeds by expanding the latter quantities in terms of the former, as
\begin{equation}
    \hom(G, \Phi(H)) = \sum_{G^{\prime} \in \mathcal{G}} c_{\Phi, G, G^{\prime}} \hom(G^{\prime}, H) \label{eq:expansion}
\end{equation}
for some coefficients $c_{\Phi, G, G^{\prime}}$, only finitely many of which are non-zero for a given $G$.
Abstractly, this is possible because the homomorphism counts span all polynomials of the adjacency matrix of a graph that are invariant under graph isomorphism (i.e., under the two-sided action of the symmetric group on matrices)~\cite{Thi00,CDM17}.
Since $\hom(G, \Phi(H))$ is such a polynomial of the adjacency matrix of $H$ for those $\Phi$ we consider, it must admit an expansion of the above form for a given dimension $|V(H)|$.
Further, for those $\Phi$ we consider, we show that the coefficients $c_{\Phi, G, G^{\prime}}$ are independent of $|V(H)|$.

In order to prove results of the above form, it suffices to show that if the treewidth of $G$ is at most $k$, then all $G^{\prime}$ for which $c_{\Phi, G, G^{\prime}} \neq 0$ must have treewidth at most $t$.
Thus, our question reduces to one of the structure of the above expansion, which we address using the tools of \emph{gadget expansion} and \emph{gadget decomposition} introduced in \cref{sec:gadgets}.
There we also show that it is possible to reason about such expansions in a convenient \emph{diagrammatic} fashion, deducing which $G^{\prime}$ appear in such an expansion by manipulating $G$ in appropriate ways.
We give such upper bound arguments in \cref{sec:upper,sec:k1} for line graphs and \cref{subsec:kge2-rs,subsec:k1-rs} for $(r, S)$-line graphs.

For lower bounds on the transfer parameter, we must show a converse result, exhibiting pairs of graphs $H_1, H_2$ that agree on homomorphism counts of graphs $G$ of treewidth at most $t$ but for which $\Phi(H_1), \Phi(H_2)$ disagree on the homomorphism count of some graph of treewidth $k$.
In the parameter ranges where our bounds are sharp, we first show that the expansion \eqref{eq:expansion} indeed contains some graph of the largest possible treewidth (below our upper bounds) with non-zero coefficients.
We then combine this observation with the theory of \emph{homomorphism-distinguishing closedness} of graphs of bounded treewidth \cite{Neuen23,Seppelt23}.
The simplest form of this property says that every graph of treewidth greater than $t$ distinguishes some pair $H_1, H_2$ that agree on homomorphism counts from all graphs of treewidth at most $t$.
A more general form for linear combinations of homomorphism counts (\cref{thm:hdc}) guarantees the same conclusion for any finite linear combination containing such a graph with non-zero coefficient.
Applying this to \eqref{eq:expansion} for a suitable $G$ of treewidth at most $k$ gives the required pair $H_1, H_2$.
We provide these arguments in \cref{sec:lower} for line graphs and \cref{subsec:kge2-rs,subsec:k1-rs,subsec:cycle} for $(r, S)$-line graphs.

The idea of taking a diagrammatic decomposition of the effect of graph operations on adjacency matrices is closely related to the constructions of \emph{graph matrices}~\cite{AMP16,BHK19} and \emph{tensor networks}~\cite{BB17,MW19} appearing in other applications in combinatorics, linear algebra, and theoretical computer science.
We believe that this perspective on graph operations may be of independent interest and could shed light on other questions about their properties.

\subsection{Related work}\label{subsec:related-work}
Homomorphism counts from all source graphs characterize graph isomorphism~\cite{Lovasz67}, and restricting the source graphs to treewidth at most $k$ characterizes $\WLp{k}$-equivalence~\cite{Dvorak10,DGR18}.
Related characterizations use bounded pathwidth and treedepth~\cite{GRS22}, or planar graphs~\cite{MR20}, and are collected in Seppelt's overview~\cite{HomIndZoo}.
Our gadget expansions use the algebra of bi-labeled graphs and homomorphism matrices~\cite{MR20,GRS22}, recalled in \cref{subsec:algebra}, to express how graph operations transform these counts. Recall that \emph{quantum graphs} are formal finite linear combinations of graphs on which $\hom(\cdot,H)$ acts linearly~\cite[Chapter~6]{Lovasz12}. In this language, \eqref{eq:expansion} states that $\hom(G,\Phi(H))=\hom(Q_{\Phi,G},H)$ for the quantum graph $Q_{\Phi,G}=\sum_{G'}c_{\Phi,G,G'}\,G'$.

Our upper bounds also follow from finite-variable counting logic. $\WLp{k}$-equivalence coincides with equivalence in the counting logic $\mathsf C^{k+1}$ with $k+1$ variables~\cite{CFI92,Dvorak10}. Otto's lemma on interpretations~\cite[Lemma~1.50]{Otto97} shows that $\mathsf C^{k+1}$-definable properties of an interpretation in the $s$-th power are $\mathsf C^{s(k+1)}$-definable properties of the input. Because the $(r,S)$-line graph is an interpretation in the $r$-th power, this yields $\tau_{r,S}(k)\le rk+r-1$, and in particular $\tau_L(k)\le2k+1$. We provide a proof via counting homomorphisms and gadget expansions.

Our lower bounds combine Neuen's homomorphism-distinguishing closedness theorem for the class of graphs of bounded treewidth~\cite{Neuen23} with Seppelt's lemma for linear combinations~\cite{Seppelt23}.
Together, these results turn a non-zero coefficient of a graph of sufficiently large treewidth into a pair witnessing failure of WL-equivalence after a graph operation.
This applies the theory initiated by Roberson~\cite{Roberson22}, who conjectured closedness for every graph class closed under minors and disjoint unions.
Related closedness results include planar graphs~\cite{Roberson22}, graphs of bounded treedepth~\cite{FSS24}, and outerplanar and cactus graphs~\cite{NS26}.

Prior work on graph learning has used representations based on edge adjacencies, including line graph neural networks~\cite{ChenLiBruna19} and dual message passing~\cite{LiuSong22}.
In the same spirit as our motivation from ``boosting'' methods for refuting graph isomorphism, Yang and Huang~\cite{YangHuang24} study line graph transformations of strongly regular graphs~\cite{GroheNeuen21} and Cai--F\"urer--Immerman graphs~\cite{CFI92}, two standard sources of graph pairs that fixed levels of WL fail to distinguish.
In recent work, Yang~\cite{Yang26} mainly studies the converse question of how WL-equivalence of line graphs transfers to WL-equivalence of the original graphs. His Proposition~2 is the case $k=1$ of \cref{thm:main-lg}.

\subsection{Notation}

All graphs are finite, simple, and undirected unless otherwise specified.
For ease of notation we denote the undirected edge $\{u,v\}$ by $uv$ or $vu$ interchangeably.
We freely identify a graph $G$ with its adjacency matrix and denote both by $G$.

\section{Preliminaries}\label{sec:prelim}

\subsection{Line graphs and Krausz cliques}\label{subsec:krausz}

Recall the line graph operation of \cref{defn:lg}.
For context, we mention the following classical theorem showing it is essentially invertible.

\begin{thm}[Whitney~\cite{Whitney32}]\label{thm:whitney}
    If $H_1$ and $H_2$ are connected graphs with $L(H_1)\cong L(H_2)$, then
    $H_1\cong H_2$, unless $\{H_1,H_2\}=\{K_3,K_{1,3}\}$.
\end{thm}

\noi Our lower bound constructions (\cref{lem:realization-lg,lem:realization-rs}) use the following Krausz cliques; for context we also recall Krausz's characterization of line graphs.

\begin{defn}[Krausz clique]\label{defn:krausz-clique}
    For $v\in V(H)$ write
    \[
        C_v=\{\,e\in E(H):v\in e\,\}\subseteq V(L(H))
    \]
    for the \emph{Krausz clique} of $v$. Any two edges
    through $v$ are adjacent in $L(H)$, so $C_v$ is indeed a clique in $L(H)$. Moreover:
    \begin{itemize}
        \item Each vertex $e=uv$ of $L(H)$ lies in exactly two Krausz cliques, namely $C_u$ and $C_v$.
        \item Each edge $ef$ of $L(H)$ lies in exactly one Krausz clique, namely $C_v$ where $v$ is the
        shared endpoint of $e$ and $f$ in $H$.
    \end{itemize}
\end{defn}

\begin{prop}[Krausz partition \cite{Krausz43}]\label{prop:Krausz-partition}
    A graph $H$ is a line graph if and only if it admits an edge partition into cliques with every vertex
    in at most two of them.
\end{prop}

\subsection{Treewidth}

Recall that treewidth appears in the definition of WL-equivalence (\cref{defn:wl}).

\begin{defn}[Tree decomposition]\label{defn:tree-decomp}
    A \emph{tree decomposition} of a graph $G$ consists of a tree $T$ and a collection of subsets $\{B_t\}_{t\in V(T)}$ of $V(G)$, called \emph{bags}, such that
    \begin{enumerate}
        \item Every vertex and every edge of $G$ is contained in some bag.
        \item For every vertex $v$ of $G$ the set of bags containing $v$ induces a connected subtree of $T$.
    \end{enumerate}
    The \emph{width} of a tree decomposition $(T, \{B_t\}_{t\in V(T)})$ is $\max_{t\in V(T)}|B_t|-1$.
    The \emph{treewidth} of $G$, denoted $\tw(G)$, is the minimum such width. Note $\tw(G)\le1$ exactly when $G$ is a forest. Write $\T_k\subseteq\G$ for the set of isomorphism classes of simple graphs with treewidth at most $k$.
\end{defn}

\subsection{Homomorphism indistinguishability}\label{subsec:hom}

For $G,H\in\G$, we write $\hom(G,H)$ for the number of \emph{homomorphisms} $V(G)\to V(H)$, i.e., maps sending edges of $G$
to edges of $H$.

\begin{defn}\label{defn:hom-indist}
    Two graphs $H_1,H_2$ are \emph{homomorphism-indistinguishable} over a class of graphs $\mathcal C$ if $\hom(G,H_1)=\hom(G,H_2)$
    for all $G\in\mathcal C$.
\end{defn}

\noindent
A classical result of Lov\'{a}sz~\cite{Lovasz67} states that isomorphism is equivalent to homomorphism-indistinguishability
over the class of all graphs.
Similarly, WL-equivalence is characterized by homomorphism-indistinguishability over $\T_k$~\cite{Dvorak10,DGR18}, as in \cref{defn:wl}.

\begin{defn}\label{defn:hdc}
    A class of graphs $\mathcal C$ is \emph{homomorphism-distinguishing closed} if for every $G\not\in\mathcal C$ there exist $H_1, H_2$
    homomorphism-indistinguishable over $\mathcal C$ with $\hom(G,H_1)\neq\hom(G,H_2)$.
\end{defn}

\noi The following theorem is critical for the lower bound of \cref{thm:main-lg}; in \cref{sec:lower} we use its
contrapositive to show the existence of $H_1\WL{2k}H_2$ with $L(H_1)\not\WL{k}L(H_2)$.

\begin{thm}[Theorem 2 of \cite{Neuen23} and Lemma 6 of \cite{Seppelt23}, arXiv v4]\label{thm:hdc}
    For every $k\ge1$ the class $\T_k$ is homomorphism-distinguishing closed.
    Consequently, if $\Gamma$ is a finite set of pairwise non-isomorphic graphs and
    $\beta\colon\Gamma\to\R\setminus\{0\}$ satisfies
    \[
        H_1\WL{k}H_2\;\Longrightarrow\;
        \sum_{\gamma\in\Gamma}\beta(\gamma)\;\hom(\gamma,H_1)
        =\sum_{\gamma\in\Gamma}\beta(\gamma)\;\hom(\gamma,H_2)
        \qquad\text{for all }H_1,H_2,
    \]
    then $\Gamma\subseteq\T_k$.
\end{thm}

\subsection{The algebra of bi-labeled graphs}\label{subsec:algebra}
Unless otherwise specified, all functions $G\rightarrow H$ below are taken to be graph homomorphisms.
The bi-labeled graph algebra and its tensor presentation have been formalized by Ma\v{n}\v{c}inska and
Roberson~\cite{MR20} and used by Grohe, Rattan, and Seppelt~\cite{GRS22}.

\begin{defn}[Bi-labeled graph]\label{defn:bl}
    A \emph{bi-labeled graph} is a triple $(G,\boldsymbol{a},\boldsymbol{b})$
    where $G\in\G$, and $\boldsymbol{a}\in V(G)^k,\boldsymbol{b}\in V(G)^\ell$
    are tuples of vertices in $G$. We refer to $\boldsymbol{a}$ as the \emph{left interface},
    and its vertices as the \emph{left vertices} of $\boldsymbol{G}$. Similarly, $\boldsymbol{b}$
    is the \emph{right interface} and its vertices are the \emph{right vertices}
    of $\boldsymbol{G}$. $G$ is the \emph{underlying graph} of $\boldsymbol{G}$.
    We use $\boldsymbol{\G}_{k,\ell}$ to denote the set of isomorphism classes of bi-labeled graphs
    with $k$ left and $\ell$ right vertices. Then $\boldsymbol{\G}=\cup_{k,\ell}\boldsymbol{\G}_{k,\ell}$ is the
    set of all bi-labeled graphs.
\end{defn}

\noi When drawing a bi-labeled graph $\boldsymbol{G}=(G,\boldsymbol{a},\boldsymbol{b})$ we put exactly $k$ vertices on the left and $\ell$ vertices on the right. If a vertex $v\in V(G)$ occurs multiple times in the concatenated interface $\boldsymbol{a}\boldsymbol{b}$ we visualize these occurrences as distinct vertices joined by a \emph{relation edge} indicating they are the same vertex. We use dashed lines for relation edges and solid lines for ordinary graph edges.

\begin{defn}[Homomorphism of bi-labeled graphs]\label{defn:homomorphism-bl}
    Given $\boldsymbol{G},\boldsymbol{H}\in\boldsymbol{\G}_{k,\ell}$, a \emph{homomorphism of bi-labeled graphs}
    $\varphi:\boldsymbol{G}\to\boldsymbol{H}$
    is a graph homomorphism sending the left and right vertices of $\boldsymbol{G}$ to the left and right
    vertices of $\boldsymbol{H}$ respectively.
\end{defn}

\tikzset{
  vtx/.style={circle, fill=black, inner sep=1.6pt},
  hedge/.style={line width=0.7pt},
  ident/.style={line width=0.7pt, dashed},
  gtitle/.style={font=\small, align=center},
}

\newcommand{\primMoneZero}{\begin{tikzpicture}[baseline=(title.base)]
    \node[vtx,label=left:$a_1$] (x1) at (0,0.55) {};
    \node (phantom) at (1.4,0.55) {};
    \node[gtitle] (title) at (0,-0.95) {$\boldsymbol{M}^{1,0}$\\ \;};
  \end{tikzpicture}}

\newcommand{\primMoneTwo}{\begin{tikzpicture}[baseline=(title.base)]
    \node[vtx,label=left:$a_1$]  (x1) at (0,0.55)   {};
    \node[vtx,label=right:$b_1$] (y1) at (1.4,1.1)  {};
    \node[vtx,label=right:$b_2$] (y2) at (1.4,0)    {};
    \draw[ident] (x1) -- (y1);
    \draw[ident] (x1) -- (y2);
    \node[gtitle] (title) at (0.7,-0.95)
      {$\boldsymbol{M}^{1,2}$\\[-1pt]\scriptsize $a_1\!=\!b_1\!=\!b_2$};
  \end{tikzpicture}}

\newcommand{\primA}{\begin{tikzpicture}[baseline=(title.base)]
    \node[vtx,label=left:$a_1$]  (x1) at (0,0.55)   {};
    \node[vtx,label=right:$b_1$] (y1) at (1.4,0.55) {};
    \draw[hedge] (x1) -- (y1);
    \node[gtitle] (title) at (0.7,-0.95)
      {$\boldsymbol{A}$\\[-1pt]\scriptsize $a_1\!\sim\!b_1$};
  \end{tikzpicture}}

\newcommand{\primS}{\begin{tikzpicture}[baseline=(title.base)]
    \node[vtx,label=left:$a_1$]  (x1) at (0,1.1)   {};
    \node[vtx,label=left:$a_2$]  (x2) at (0,0)     {};
    \node[vtx,label=right:$b_1$] (y1) at (1.4,1.1) {};
    \node[vtx,label=right:$b_2$] (y2) at (1.4,0)   {};
    \draw[ident] (x1) -- (y2);
    \draw[ident] (x2) -- (y1);
    \node[gtitle] (title) at (0.7,-0.95)
      {$\boldsymbol{S}$\\[-1pt]\scriptsize $a_1\!=\!b_2,\ a_2\!=\!b_1$};
  \end{tikzpicture}}

\begin{figure}
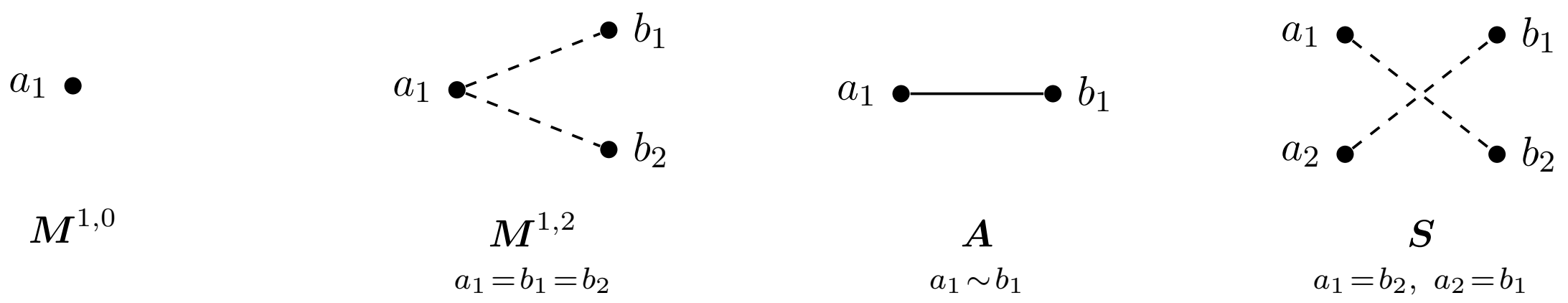

  \centering
  \makebox[\textwidth]{\hspace{0.06\textwidth}\primMoneZero
    \hfill
    \primMoneTwo
    \hfill
    \primA
    \hfill
    \primS
    \hspace{0.06\textwidth}}
  \caption{We illustrate some important bi-labeled graphs, along with the relations among vertices that they imply.}
  \label{fig:primitives}
\end{figure}

\begin{defn}[Homomorphism matrix]
    Given a bi-labeled graph $\boldsymbol{G}=(G,\boldsymbol{a},\boldsymbol{b})\in\boldsymbol{\G}_{k,\ell}$
    and a graph $H$, the \emph{homomorphism matrix}
    $T^{\boldsymbol{G}\to H}$ is the $V(H)^k\times V(H)^\ell$ matrix defined entrywise by
    \[
        T^{\boldsymbol{G}\to H}_{\boldsymbol{u},\boldsymbol{v}}= \#\{\varphi\colon G\to H : \varphi(\boldsymbol a)=\boldsymbol u,\ \varphi(\boldsymbol b)=\boldsymbol v\}.
    \]
    Concretely, letting $V(H)=[n]$ we have
    \[
        T^{\boldsymbol{G}\to H}_{\boldsymbol{u},\boldsymbol{v}}
        =\sum_{\substack{\varphi:V(G)\to[n] \\ \varphi(\boldsymbol{a})=\boldsymbol{u}
        \\ \varphi(\boldsymbol{b})=\boldsymbol{v}}}\;\prod_{xy\in E(G)}H_{\varphi(x),\varphi(y)}.
    \]
\end{defn}

\noindent Note this is similar to the definition of \emph{graph matrices} in the theoretical computer science literature \cite{AMP16,BHK19}; here the input is a binary matrix and the sum includes non-injective maps.
As with the graph matrix operation, the map $H\mapsto T^{\boldsymbol{G}\to H}$ uses the bi-labeled ``template" graph $\boldsymbol{G}$ to construct a new, possibly larger, matrix whose entries are polynomials in the entries of $H$.
Conveniently, natural operations on the bi-labeled graphs $\boldsymbol{G}$ correspond to natural linear-algebraic operations on the homomorphism matrices $T^{\boldsymbol{G}\to H}$.

\begin{defn}[Bi-labeled graph operations]\label{defn:operations-bl}
    Consider the bi-labeled graphs $\boldsymbol{G}_1=(G_1,\boldsymbol{a}_1,\boldsymbol{b}_1)\in\boldsymbol{\G}_{k_1,\ell_1}$ and
    $\boldsymbol{G}_2=(G_2,\boldsymbol{a}_2,\boldsymbol{b}_2)\in\boldsymbol{\G}_{k_2,\ell_2}$. We record the primary operations
    of the bi-labeled graph algebra below and their corresponding operations on homomorphism matrices (see \cref{fig:operations}).
    Throughout, vertex identifications retain loops and collapse parallel edges.
    \begin{itemize}
        \item \emph{\textbf{Series Composition:}}
            $\boldsymbol{G}_1\circ\boldsymbol{G}_2$ is defined when $\ell_1=k_2$ and has underlying graph formed by taking the disjoint union $G_1\sqcup G_2$ and identifying $\boldsymbol{b}_1=\boldsymbol{a}_2$. It has left interface $\boldsymbol{a}_1$ and right interface $\boldsymbol{b}_2$. We have
            \[
                T^{\boldsymbol{G}_1\circ\boldsymbol{G}_2\to H}=T^{\boldsymbol{G}_1\to H}T^{\boldsymbol{G}_2\to H},
            \]
            where juxtaposition on the right denotes the ordinary product of matrices.
        \item \emph{\textbf{Tensor Product:}}
            $\boldsymbol{G}_1\otimes\boldsymbol{G}_2$ has underlying graph $G_1\sqcup G_2$. It has left interface $\boldsymbol{a}_1\boldsymbol{a}_2$ and right interface $\boldsymbol{b}_1\boldsymbol{b}_2$ where juxtaposition denotes concatenation. We have
            \[
                T^{\boldsymbol{G}_1\otimes\boldsymbol{G}_2\to H}=T^{\boldsymbol{G}_1\to H}\otimes T^{\boldsymbol{G}_2\to H}.
            \]
        \item \emph{\textbf{Transposition:}}
            $\boldsymbol{G}_1^\top$ is formed by swapping the left and right interfaces of $\boldsymbol{G}_1$. We have
            \[
                T^{\boldsymbol{G}_1^\top\to H}=\left(T^{\boldsymbol{G}_1\to H}\right)^\top.
            \]
        \item \emph{\textbf{Parallel Composition:}}
            $\boldsymbol{G}_1\bullet\boldsymbol{G}_2$ is defined when $k_1=k_2$ and $\ell_1=\ell_2$ and has underlying graph formed by taking the disjoint union $G_1\sqcup G_2$ and identifying $\boldsymbol{a}_1=\boldsymbol{a}_2$ and $\boldsymbol{b}_1=\boldsymbol{b}_2$. It has left interface $\boldsymbol{a}_1=\boldsymbol{a}_2$ and right interface $\boldsymbol{b}_1=\boldsymbol{b}_2$. We have
            \[
                T^{\boldsymbol{G}_1\bullet\boldsymbol{G}_2\to H}=T^{\boldsymbol{G}_1\to H}\bullet T^{\boldsymbol{G}_2\to H}
            \]
            where $\bullet$ on the right denotes the entrywise or Hadamard product of matrices.
    \end{itemize}
\end{defn}

\tikzset{
  opv/.style={circle, fill=black, inner sep=1.5pt},
  opbox/.style={draw=black, line width=0.8pt, rounded corners=2pt, fill=black!4},
  opident/.style={line width=0.7pt, dashed},
  opwire/.style={line width=0.7pt},
  oplbl/.style={font=\scriptsize},
  opname/.style={font=\small\bfseries},
}

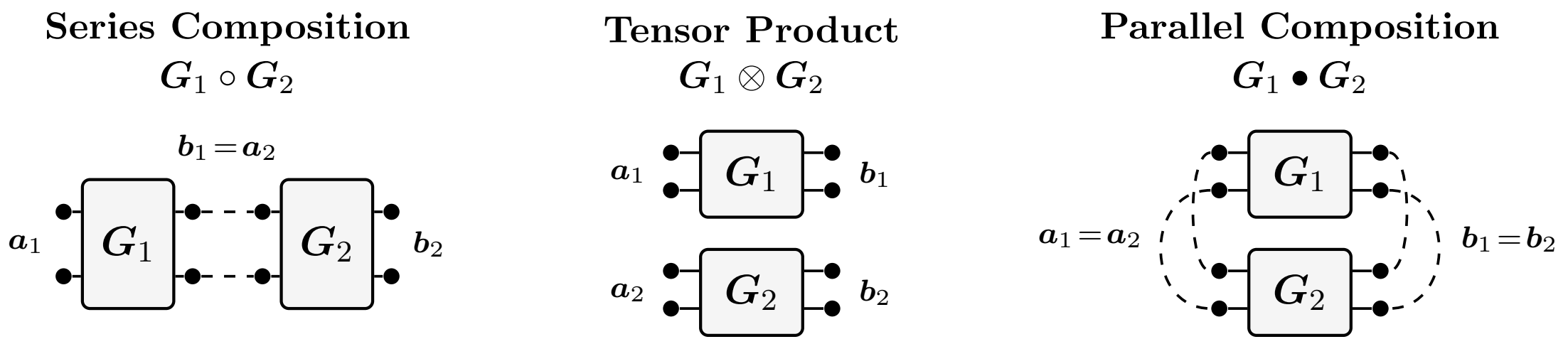
\begin{figure}
\centering
\begin{tikzpicture}

\begin{scope}[shift={(0,0)}]
  \node[opname] at (1.875,1.5) {Series Composition};
  \node[opname] at (1.875,1.05) {$\boldsymbol{G}_1\circ\boldsymbol{G}_2$};
  \node[opbox, minimum width=0.85cm, minimum height=1.2cm] (g1) at (0.95,-0.5) {$\boldsymbol{G}_1$};
  \node[opbox, minimum width=0.85cm, minimum height=1.2cm] (g2) at (2.80,-0.5) {$\boldsymbol{G}_2$};
  \node[opv] (a1u) at (0.35,-0.8) {}; \node[opv] (a1d) at (0.35,-0.2) {};
  \node[opv] (b1u) at (1.55,-0.8) {}; \node[opv] (b1d) at (1.55,-0.2) {};
  \draw[opwire] (a1u)--(g1.west|-a1u); \draw[opwire] (a1d)--(g1.west|-a1d);
  \draw[opwire] (b1u)--(g1.east|-b1u); \draw[opwire] (b1d)--(g1.east|-b1d);
  \node[opv] (a2u) at (2.20,-0.8) {}; \node[opv] (a2d) at (2.20,-0.2) {};
  \node[opv] (b2u) at (3.40,-0.8) {}; \node[opv] (b2d) at (3.40,-0.2) {};
  \draw[opwire] (a2u)--(g2.west|-a2u); \draw[opwire] (a2d)--(g2.west|-a2d);
  \draw[opwire] (b2u)--(g2.east|-b2u); \draw[opwire] (b2d)--(g2.east|-b2d);
  \draw[opident] (b1u)--(a2u); \draw[opident] (b1d)--(a2d);
  \node[oplbl] at (1.875,0.4) {$\boldsymbol{b}_1\!=\!\boldsymbol{a}_2$};
  \node[oplbl] at (0,-0.5) {$\boldsymbol{a}_1$};
  \node[oplbl] at (3.75,-0.5) {$\boldsymbol{b}_2$};
\end{scope}

\begin{scope}[shift={(10.25,0)}]
  \node[opname] at (1.6,1.5) {Parallel Composition};
  \node[opname] at (1.6,1.05) {$\boldsymbol{G}_1\bullet\boldsymbol{G}_2$};
  \node[opbox, minimum width=0.95cm, minimum height=0.8cm] (h1) at (1.6,0.15) {$\boldsymbol{G}_1$};
  \node[opbox, minimum width=0.95cm, minimum height=0.8cm] (h2) at (1.6,-0.95) {$\boldsymbol{G}_2$};
  \node[opv] (a1u) at (0.85,0.0) {}; \node[opv] (a1d) at (0.85,0.35) {};
  \node[opv] (b1u) at (2.35,0.0) {}; \node[opv] (b1d) at (2.35,0.35) {};
  \draw[opwire] (a1u)--(h1.west|-a1u); \draw[opwire] (a1d)--(h1.west|-a1d);
  \draw[opwire] (b1u)--(h1.east|-b1u); \draw[opwire] (b1d)--(h1.east|-b1d);
  \node[opv] (a2u) at (0.85,-1.1) {}; \node[opv] (a2d) at (0.85,-0.75) {};
  \node[opv] (b2u) at (2.35,-1.1) {}; \node[opv] (b2d) at (2.35,-0.75) {};
  \draw[opwire] (a2u)--(h2.west|-a2u); \draw[opwire] (a2d)--(h2.west|-a2d);
  \draw[opwire] (b2u)--(h2.east|-b2u); \draw[opwire] (b2d)--(h2.east|-b2d);
  \draw[opident] (a1u) .. controls (0.15,0.0) and (0.15,-1.1) .. (a2u);
  \draw[opident] (a1d) .. controls (0.55,0.35) and (0.55,-0.75) .. (a2d);
  \draw[opident] (b1u) .. controls (3.05,0.0) and (3.05,-1.1) .. (b2u);
  \draw[opident] (b1d) .. controls (2.65,0.35) and (2.65,-0.75) .. (b2d);
  \node[oplbl] at (-0.35,-0.45) {$\boldsymbol{a}_1\!=\!\boldsymbol{a}_2$};
  \node[oplbl] at (3.55,-0.45) {$\boldsymbol{b}_1\!=\!\boldsymbol{b}_2$};
\end{scope}

\begin{scope}[shift={(5.25,0)}]
  \node[opname] at (1.5,1.5) {Tensor Product};
  \node[opname] at (1.5,1.05) {$\boldsymbol{G}_1\otimes\boldsymbol{G}_2$};
  \node[opbox, minimum width=0.95cm, minimum height=0.8cm] (t1) at (1.5,0.15) {$\boldsymbol{G}_1$};
  \node[opbox, minimum width=0.95cm, minimum height=0.8cm] (t2) at (1.5,-0.95) {$\boldsymbol{G}_2$};
  \node[opv] (i1u) at (0.75,0.0) {}; \node[opv] (i1d) at (0.75,0.35) {};
  \node[opv] (o1u) at (2.25,0.0) {}; \node[opv] (o1d) at (2.25,0.35) {};
  \draw[opwire] (i1u)--(t1.west|-i1u); \draw[opwire] (i1d)--(t1.west|-i1d);
  \draw[opwire] (o1u)--(t1.east|-o1u); \draw[opwire] (o1d)--(t1.east|-o1d);
  \node[opv] (i2u) at (0.75,-1.1) {}; \node[opv] (i2d) at (0.75,-0.75) {};
  \node[opv] (o2u) at (2.25,-1.1) {}; \node[opv] (o2d) at (2.25,-0.75) {};
  \draw[opwire] (i2u)--(t2.west|-i2u); \draw[opwire] (i2d)--(t2.west|-i2d);
  \draw[opwire] (o2u)--(t2.east|-o2u); \draw[opwire] (o2d)--(t2.east|-o2d);
  \node[oplbl] at (0.35,0.15) {$\boldsymbol{a}_1$};
  \node[oplbl] at (0.35,-0.95) {$\boldsymbol{a}_2$};
  \node[oplbl] at (2.65,0.15) {$\boldsymbol{b}_1$};
  \node[oplbl] at (2.65,-0.95) {$\boldsymbol{b}_2$};
\end{scope}

\end{tikzpicture}
\caption{We illustrate the effect of the three binary operations described in \cref{defn:operations-bl} on pairs of bi-labeled graphs.}
\label{fig:operations}
\end{figure}
\begin{thm}[Theorem 8.4 of \cite{MR20}]
    $\boldsymbol{\G}=\langle \boldsymbol{M}^{1,0},\boldsymbol{M}^{1,2},\boldsymbol{S},\boldsymbol{A}\rangle_{\circ,\otimes,\top}$. That is, all bi-labeled graphs can be generated from the four bi-labeled graphs given in \cref{fig:primitives} by combining them with series composition, tensor product, and transposition.
\end{thm}
\noindent We do not use this result, but it shows that the gadgets below are built from standard primitives.

\section{Gadget expansions and decompositions}\label{sec:gadgets}

We now present the tools, based on the algebra of bi-labeled graphs, used for realizing graph operations via diagrammatic expansions.

\begin{defn}[Gadget, gadget family, interface size]\label{defn:gadget}
    A \emph{gadget family} $(\F,\alpha)$ is a finite set of bi-labeled graphs $\F$, each with left and right interfaces of size $\mu$, together with a \emph{coefficient function} $\alpha:\F\to\Z\setminus\{0\}$. We refer to a bi-labeled graph $\boldsymbol{F}\in\F$ as a \emph{gadget} and $\mu$ as the \emph{interface size} of the family $\F$.
\end{defn}

\begin{defn}[Transpose-closed]
    We say a gadget family $(\F,\alpha)$ is \emph{transpose-closed} when
    \[
    \boldsymbol{F}\in\F\;\Longrightarrow\;\boldsymbol{F}^\top\in\F,\;\alpha(\boldsymbol{F}^\top)=\alpha(\boldsymbol{F}).
    \]
\end{defn}

Given a graph $G$, the following definition describes how we can derive a new graph by ``plugging in'' a gadget along each edge of $G$.

\begin{defn}[Gadget expansion]\label{defn:gex}
    Let $\F$ be a gadget family of interface size $\mu$, $G$ be a directed graph, and $g\colon E(G)\to\F$ be a \emph{gadget
    assignment}. The \emph{gadget expansion} $\gex_g(G)$ is a graph, possibly with loops, constructed by inserting a copy of $g(e)$ along each edge $e\in E(G)$. Formally, for each vertex $u\in V(G)$ create an interface $\boldsymbol{u}=(u_1,\hdots,u_\mu)$. Then for each directed edge $e=uv\in E(G)$ create a disjoint copy of $g(e)$; identify the left and right interfaces of $g(e)$ with $\boldsymbol{u}$ and $\boldsymbol{v}$ respectively; and collapse any resulting parallel edges, retaining loops. The result is
    \[
        V(\gex_g(G))=\left[\left(V(G)\times[\mu]\right)\sqcup\bigsqcup_e V(\boldsymbol{F}_e)\right]/\sim\,,\qquad E(\gex_g(G))=\left[\bigsqcup_eE(\boldsymbol{F}_e)\right]/\sim
    \]
    where $\sim$ is the equivalence relation generated by the interface identifications.
\end{defn}

\begin{lem}[Orientation independence]\label{lem:orientation-indep}
    Let $(\F,\alpha)$ be transpose-closed, $G$ be an undirected graph, and
    $o,o'$ be two orientations of $G$. Given an assignment $g\colon E(G)\to\F$ under
    $o$, define $g'$ under $o'$ by $g'(e)=g(e)^\top$ on edges where $o,o'$ disagree and $g'(e)=g(e)$ elsewhere. Then $g\mapsto g'$ is a bijection between
    assignments under $o$ and $o'$, and letting $G_o$ denote the orientation of $G$ by $o$,
    \[
        \gex_{g}(G_{o})=\gex_{g'}(G_{o'}),\qquad
        \prod_{e\in E(G)}\alpha\left(g(e)\right)=\prod_{e\in E(G)}\alpha\left(g'(e)\right).
    \]
\end{lem}
\noindent
Thus changing the orientation of $G$ merely reindexes the gadget assignments, preserving each expanded graph and its coefficient.
In particular, the weighted sum over all gadget assignments is independent of the chosen orientation.
Henceforth, gadget families are assumed to be transpose-closed, and we choose an arbitrary orientation whenever applying the construction to an undirected graph.

\begin{proof}
    Every pair of orientations $o,o'$ is connected by a sequence of single-edge
    reversals, so it suffices to consider the case where $o,o'$ differ on a single
    edge. Let $e=\{u,v\}$ with $o(e)=uv$ and $o'(e)=vu$, and write
    $g(e)=(F,\boldsymbol a,\boldsymbol b)$, so $g'(e)=g(e)^\top=(F,\boldsymbol
    b,\boldsymbol a)$. Under $o$, the expansion $\gex_g(G_o)$ identifies the left
    interface $\boldsymbol a$ of $g(e)$ with $\boldsymbol u$ and the right interface
    $\boldsymbol b$ with $\boldsymbol v$. Under $o'$ the edge points $v\to u$, so
    $\gex_{g'}(G_{o'})$ identifies the left interface of $g'(e)$ with $\boldsymbol v$
    and its right interface with $\boldsymbol u$ since $g'(e)$ has left interface
    $\boldsymbol b$ and right interface $\boldsymbol a$. In both cases the identifications are the same, and the underlying
    graph $F$ is unchanged, so $\gex_g(G_o)=\gex_{g'}(G_{o'})$. Given an assignment $g'$ under $o'$, the same construction yields the inverse map, hence $g\mapsto g'$ is a bijection. Finally, transpose closure gives $\alpha(g(e)^\top)=\alpha(g(e))$, so the products coincide.
\end{proof}

We now describe how to compute homomorphism counts of graphs formed as gadget expansions in terms of the underlying gadgets.

\begin{defn}[Contraction along $G$]\label{defn:contraction}
    Let $G$ be a directed graph and fix integers $n,\mu\ge1$. Given an edge-indexed
    tuple $\mathcal A=(A_e)_{e\in E(G)}$ of $[n]^\mu\times[n]^\mu$ matrices, the
    \emph{contraction of $\mathcal A$ along $G$} is
    \[
        Z_G(\mathcal A)
        =\sum_{\varphi:V(G)\to[n]^\mu}\ \prod_{e=uv\in E(G)}(A_e)_{\varphi(u),\varphi(v)}.
    \]
    The map $\mathcal A\mapsto Z_G(\mathcal A)$ is multilinear in $(A_e)_{e\in E(G)}$.
    When $\mu=1$ and $A_e\equiv H$ for $H$ the adjacency matrix of a graph, we recover
    $Z_G(\mathcal A)=\hom(G,H)$.
    For an undirected graph $G$ and a symmetric $[n]^\mu\times[n]^\mu$ matrix $A$, we define $\hom(G,A)=Z_G\bigl((A)_{e\in E(G)}\bigr)$ for any orientation of $G$; by symmetry of $A$ this does not depend on the orientation.
\end{defn}

\begin{lem}[Gadget hom count]\label{lem:gex-count}
    Let $G$ be a directed graph, $H$ a simple undirected graph on $[n]$, $\F$ a gadget family of interface size $\mu$, and $g\colon E(G)\to\F$ an assignment. Then
    \[
        \hom(\gex_g(G),H)
        =\sum_{\varphi:V(G)\to [n]^\mu}\ \prod_{e=uv\in E(G)} T^{g(e)\to H}_{\varphi(u),\varphi(v)}=Z_G\left(\bigl(T^{g(e)\to H}\bigr)_{e}\right).
    \]
\end{lem}

\begin{proof}
    The second equality is the definition of $Z_G$. To prove the first equality, first identify each
    map $\varphi\colon V(G)\to V(H)^\mu$ with the map $\psi\colon V(G)\times[\mu]\to V(H)$
    given by $\varphi(u)=(\psi(u_1),\dots,\psi(u_\mu))$. Fix an edge $e=uv$ and write $g(e)=(F_e,\boldsymbol a_e,\boldsymbol b_e)$. By the
    definition of the homomorphism matrix,
    \begin{equation}
        T^{g(e)\to H}_{\varphi(u),\varphi(v)}
        =\sum_{\theta_e\in \Theta_e(\psi)}\ \prod_{xy\in E(F_e)} H_{\theta_e(x),\theta_e(y)},
    \end{equation}
    where $\Theta_e(\psi)$ is the set of maps $\theta_e\colon V(F_e)\to V(H)$ with
    $\theta_e(\boldsymbol a_e)=\varphi(u)$ and $\theta_e(\boldsymbol b_e)=\varphi(v)$. Then
    \[
        \sum_{\varphi}\;\prod_{e\in E(G)}
            T^{g(e)\to H}_{\varphi(u),\varphi(v)}
        =\sum_{\psi}\;\prod_{e\in E(G)}\;
            \sum_{\theta_e\in \Theta_e(\psi)}\ \prod_{xy\in E(F_e)} H_{\theta_e(x),\theta_e(y)}.
    \]
    Distinct underlying graphs $F_e$ share only interface vertices, whose values are fixed by
    $\psi$. Thus the free variables of the terms of the inner sum are disjoint, and distributivity yields
    \[
        \prod_{e\in E(G)}\;\sum_{\theta_e\in \Theta_e(\psi)}\ \prod_{xy\in E(F_e)} H_{\theta_e(x),\theta_e(y)}
        =\sum_{(\theta_e)_e\,\in\,\prod_e \Theta_e(\psi)}\;
        \prod_{e\in E(G)}\ \prod_{xy\in E(F_e)} H_{\theta_e(x),\theta_e(y)}.
    \]
    A choice of $\psi$ together with a family $(\theta_e)_{e\in E(G)}$ assigns a vertex of
    $H$ to every element of
    \[
        \big(V(G)\times[\mu]\big)\sqcup\left(\bigsqcup_{e}V(F_e)\right).
    \]
    Now let $\sim$ denote the equivalence relation generated by the interface identifications of $\gex_g(G)$. By \cref{defn:gex},
    \[
        V(\gex_g(G))=\left[\left(V(G)\times[\mu]\right)\sqcup\bigsqcup_e V(F_e)\right]/\sim.
    \]
    The condition $\theta_e\in \Theta_e(\psi)$ means $\theta_e$ is consistent with the interface identifications of $\gex_g(G)$, i.e. constant on equivalence classes of $\sim$. The pairs $(\psi,(\theta_e)_{e\in E(G)})$ are thus in bijection with maps
    $\eta\colon V(\gex_g(G))\to V(H)$.
    A gadget edge $xy\in E(F_e)$ and its image in $\gex_g(G)$ receive the same
    endpoints under $\eta$, and every edge of $\gex_g(G)$ is such an image. Since $H$
    is simple every factor lies in $\{0,1\}$, so both sides of
    \[
        \prod_{e\in E(G)}\ \prod_{xy\in E(F_e)} H_{\theta_e(x),\theta_e(y)}
        =\prod_{xy\in E(\gex_g(G))} H_{\eta(x),\eta(y)}
    \]
    are indicators that $\eta$ maps every gadget edge into
    $E(H)$ and therefore agree. Summing over all $\eta\colon V(\gex_g(G))\to V(H)$ gives
    $\hom(\gex_g(G),H)$.
\end{proof}

Finally, we arrive at the main point of these tools. The definition below makes precise the sense in which a gadget family can realize a graph operation in terms of homomorphism counts. This realization allows us to express homomorphism counts from $G$ into a derived graph $\Phi(H)$ in terms of homomorphism counts from gadget expansions $\gex_g(G)$ into the original graph $H$.

\begin{defn}[Representation of graph operation, duplication factor]\label{defn:represents}
    A family $\F$ \emph{represents} a graph operation $\Phi$ with
    \emph{duplication factor} $d$ if, for every graph $H$ and every graph $G$ without isolated
    vertices,
    \begin{equation}\label{eq:represents}
        \hom\left(G,\Phi(H)\right)=d^{\,-|V(G)|}\hom\left(G,\sum_{\boldsymbol{F}\in\F}\alpha(\boldsymbol{F})\;T^{\boldsymbol{F}\to H}\right).
    \end{equation}
    Here $\sum_{\boldsymbol F\in\F}\alpha(\boldsymbol F)\,T^{\boldsymbol F\to H}$ is symmetric because $\F$ is transpose-closed.
\end{defn}

\begin{lem}[Gadget hom decomposition]\label{lem:gadget-decomp}
    If $\F$ has interface size $\mu$ and represents $\Phi$ with duplication factor $d$, then for every graph $G$ without
    isolated vertices and $H$ on $[n]$,
    \begin{equation}\label{eqn:decomp}
        \hom\left(G,\Phi(H)\right)
        =d^{-|V(G)|}\sum_{g\colon E(G)\to\F}
            \left(\prod_{e\in E(G)}\alpha\left(g(e)\right)\right)
            \hom\left(\gex_g(G),H\right).
    \end{equation}
    Letting $\Gamma\subseteq\G$ denote the set of isomorphism classes of loopless expansions $\gex_g(G)$, we have
    \begin{equation}\label{eqn:decomp-iso}
        \hom\left(G,\Phi(H)\right)
        =d^{-|V(G)|}\sum_{\gamma\in\Gamma}\beta(G,\gamma)\,\hom(\gamma,H),
    \end{equation}
    with coefficients defined for $\gamma\in\G$ by
    \begin{equation}\label{eqn:coeff-iso}
        \beta(G,\gamma)=\sum_{g\,:\,\gex_g(G)\cong\gamma}\;
            \prod_{e\in E(G)}\alpha(g(e))\;\in\;\Z.
    \end{equation}
\end{lem}

\begin{proof}
    By~\cref{defn:represents},
    \[
        d^{\,|V(G)|}\hom\left(G,\Phi(H)\right)
        =\hom\left(G,\sum_{\boldsymbol{F}\in\F}\alpha(\boldsymbol{F})\,T^{\boldsymbol{F}\to H}\right)
        =Z_G\left(\bigl(\sum_{\boldsymbol{F}\in\F}\alpha(\boldsymbol{F})\,T^{\boldsymbol{F}\to H}\bigr)_{e}\right),
    \]
    where the argument of $Z_G$ is a constant tuple. Expanding by the multilinearity of $Z_G$ in
    the tuple (\cref{defn:contraction}),
    \[
        d^{\,|V(G)|}\hom\left(G,\Phi(H)\right)
        =\sum_{g\colon E(G)\to\F}\left(\prod_{e\in E(G)}\alpha\left(g(e)\right)\right)
            \,Z_G\left(\bigl(T^{g(e)\to H}\bigr)_{e}\right)
    \]
    \[
        =\sum_{g\colon E(G)\to\F}\left(\prod_{e\in E(G)}\alpha\left(g(e)\right)\right)
            \,\hom\left(\gex_g(G),H\right),
    \]
    the last equality by \cref{lem:gex-count}.
    Dividing by $d^{\,|V(G)|}$ gives \eqref{eqn:decomp}.
    Expansions containing loops contribute zero since $H$ is simple.
    Omitting these terms and regrouping the rest by isomorphism type gives
    \eqref{eqn:decomp-iso}, with $\beta(G,\gamma)$ as in \eqref{eqn:coeff-iso}.
\end{proof}

\section{Line graphs: Proof of \texorpdfstring{\cref{thm:main-lg}}{the line-graph transfer theorem}}
\subsection{Line graph gadget family}\label{subsec:gadgets-lg}

We first describe how to apply the gadget family techniques of \cref{sec:gadgets} to the line graph operation.

\begin{defn}[Ordered line graph]\label{defn:ordered-lg}
    Given a graph $H$, the \emph{ordered line graph} $\widetilde L(H)$ is an undirected graph with vertex set $V(H)\times V(H)$ and adjacency $(a_1,a_2)\sim(b_1,b_2)$ when $a_1a_2,b_1b_2\in E(H)$ and $|\{a_1,a_2\}\cap\{b_1,b_2\}|=1$.
    Equivalently,
    \[
        \widetilde L(H)_{(a_1,a_2),(b_1,b_2)}
        =H_{a_1a_2}H_{b_1b_2}\bone\bigl[|\{a_1,a_2\}\cap\{b_1,b_2\}|=1\bigr].
    \]
\end{defn}

\begin{defn}[Line graph gadget family]\label{defn:family-lg}
    Denote the \emph{line graph gadget family} (of interface size~$2$) by $\LL$, and let $\boldsymbol{F}=(F,(a_1,a_2),(b_1,b_2))$ for $\boldsymbol{F}\in\LL$. Then $\LL$ is shown in \cref{fig:gadgets-lg} and its coefficients are given by
    \begin{center}
    \begin{tabular}{c|cccccc}
        $\boldsymbol{F}$ & $\boldsymbol{F}_{11}$ & $\boldsymbol{F}_{12}$ & $\boldsymbol{F}_{21}$ & $\boldsymbol{F}_{22}$ & $\boldsymbol{F}_\parallel$ & $\boldsymbol{F}_\times$ \\
        \hline
        $\alpha$ & $+1$ & $+1$ & $+1$ & $+1$ & $-2$ & $-2$
    \end{tabular}
    \end{center}
    where $E(F)=\{a_1a_2,b_1b_2\}$ for all $\boldsymbol{F}\in\LL$; the four \emph{single gadgets} $\boldsymbol{F}_{ij}$ identify $a_i=b_j$; and the two \emph{double gadgets}
    $\boldsymbol{F}_\parallel$ and $\boldsymbol{F}_\times$ identify $a_1=b_1$, $a_2=b_2$ and $a_1=b_2$, $a_2=b_1$ respectively.
\end{defn}

\tikzset{
  vtx/.style={circle, fill=black, inner sep=1.6pt},
  hedge/.style={line width=0.7pt},
  ident/.style={line width=0.7pt, dashed},
  gtitle/.style={font=\small, align=center},
}
\newcommand{\preframe}{\node[vtx,label=left:$a_1$]  (x1) at (0,1.1)   {};
  \node[vtx,label=left:$a_2$]  (x2) at (0,0)     {};
  \node[vtx,label=right:$b_1$] (y1) at (1.4,1.1) {};
  \node[vtx,label=right:$b_2$] (y2) at (1.4,0)   {};
  \draw[hedge] (x1)--(x2);
  \draw[hedge] (y1)--(y2);
}
\newcommand{\pregadget}[2]{\begin{tikzpicture}[baseline=(title.base)]
    \preframe #1
    \node[gtitle] (title) at (0.7,-0.95) {#2};
  \end{tikzpicture}}

\begin{figure}
  \centering
  \begin{tabular}{@{}c@{\hspace{0.1\textwidth}}c@{\hspace{0.1\textwidth}}c@{}}
    \pregadget{\draw[ident](x1)--(y1);}{$\boldsymbol{F}_{11}$,\ $\alpha=+1$\\[-1pt]\scriptsize $a_1\!=\!b_1$} &
    \pregadget{\draw[ident](x2)--(y1);}{$\boldsymbol{F}_{21}$,\ $\alpha=+1$\\[-1pt]\scriptsize $a_2\!=\!b_1$} &
    \pregadget{\draw[ident](x1)--(y1);\draw[ident](x2)--(y2);}{$\boldsymbol{F}_{\parallel}$,\ $\alpha=-2$\\[-1pt]\scriptsize $a_1\!=\!b_1,\ a_2\!=\!b_2$} \\[2em]
    \pregadget{\draw[ident](x1)--(y2);}{$\boldsymbol{F}_{12}$,\ $\alpha=+1$\\[-1pt]\scriptsize $a_1\!=\!b_2$} &
    \pregadget{\draw[ident](x2)--(y2);}{$\boldsymbol{F}_{22}$,\ $\alpha=+1$\\[-1pt]\scriptsize $a_2\!=\!b_2$} &
    \pregadget{\draw[ident](x1)--(y2);\draw[ident](x2)--(y1);}{$\boldsymbol{F}_{\times}$,\ $\alpha=-2$\\[-1pt]\scriptsize $a_1\!=\!b_2,\ a_2\!=\!b_1$} \\
  \end{tabular}
  \caption{The line graph gadget family of \cref{defn:family-lg}.}
  \label{fig:gadgets-lg}
\end{figure}

\begin{prop}\label{prop:ie-lg}
    For every graph $H$,
    \begin{equation}\label{eq:ie-lg}
        \widetilde L(H)=\sum_{\boldsymbol{F}\in\LL}\alpha(\boldsymbol{F})\;T^{\boldsymbol{F}\to H}.
    \end{equation}
\end{prop}

\begin{proof}
    Let $B_H$ and $J_s$, for $s\in\{1,2\}$, have rows and columns indexed by $V(H)^2$ and entries
    \begin{align*}
        (B_H)_{(a_1,a_2),(b_1,b_2)} &\colonequals H_{a_1a_2}H_{b_1b_2},\\
        (J_s)_{(a_1,a_2),(b_1,b_2)} &\colonequals \bone\bigl[|\{a_1,a_2\}\cap\{b_1,b_2\}|=s\bigr].
    \end{align*}
    Then $\widetilde L(H)=B_H\bullet J_1$, where $\bullet$ denotes the entrywise product.
    For $i,j\in\{1,2\}$, define $D_{ij}$ on the same index set by $(D_{ij})_{(a_1,a_2),(b_1,b_2)}\colonequals\bone[a_i=b_j]$.
    Since $H$ is loopless, every entry where $B_H$ is nonzero has $a_1\ne a_2$ and $b_1\ne b_2$.
    On these entries, $\sum_{i,j=1}^2D_{ij}$ counts the common vertices, and an intersection of size two corresponds to exactly one of the two matchings of the endpoints.
    Thus
    \begin{align*}
        B_H\bullet\sum_{i,j=1}^2D_{ij} &= B_H\bullet(J_1+2J_2),\\
        B_H\bullet J_2 &= B_H\bullet(D_{11}\bullet D_{22}+D_{12}\bullet D_{21}).
    \end{align*}
    It follows that
    \[
        \widetilde L(H)=B_H\bullet J_1
        =B_H\bullet\left(\sum_{i,j=1}^2D_{ij}-2D_{11}\bullet D_{22}-2D_{12}\bullet D_{21}\right).
    \]
    Now $B_H=T^{\boldsymbol F_0\to H}$ for the bi-labeled graph $\boldsymbol F_0=(F_0,(a_1,a_2),(b_1,b_2))$ consisting of the two disjoint edges $a_1a_2$ and $b_1b_2$.
    Entrywise multiplication by $D_{ij}$ is parallel composition with a bi-labeled graph that identifies the interface vertices $a_i$ and $b_j$.
    The six terms in the expansion above therefore correspond to the gadgets in $\LL$ as follows:
    \begin{center}
    \begin{tabular}{c|c|c|c}
        term & identification & gadget & $\alpha$\\
        \hline
        $D_{11}$ & $a_1{=}b_1$ & $\boldsymbol F_{11}$ & $+1$\\
        $D_{12}$ & $a_1{=}b_2$ & $\boldsymbol F_{12}$ & $+1$\\
        $D_{21}$ & $a_2{=}b_1$ & $\boldsymbol F_{21}$ & $+1$\\
        $D_{22}$ & $a_2{=}b_2$ & $\boldsymbol F_{22}$ & $+1$\\
        $D_{11}{\bullet} D_{22}$ & $a_1{=}b_1,a_2{=}b_2$ & $\boldsymbol F_\parallel$ & $-2$\\
        $D_{12}{\bullet} D_{21}$ & $a_1{=}b_2,a_2{=}b_1$ & $\boldsymbol F_\times$ & $-2$\\
    \end{tabular}
    \end{center}
    Summing these six terms proves the identity.
\end{proof}

\begin{lem}[Duplication factor]\label{lem:duplication-lg}
    For every graph $H$ and every graph $G$ without isolated vertices,
    \[
        \hom\left(G,\widetilde L(H)\right)=2^{|V(G)|}\,\hom\left(G,L(H)\right),
    \]
    so $\LL$ represents $L$ with duplication factor $d=2$.
\end{lem}
\begin{proof}
    Each unordered vertex $\{a,b\}$ in $L(H)$ has two corresponding ordered vertices
    $ab,ba$ in $\widetilde L(H)$, and adjacency in $\widetilde L(H)$ is the unordered
    condition $|\{a,b\}\cap\{c,d\}|=1$. Thus $\widetilde L(H)$ is obtained from $L(H)$ by
    replacing each unordered vertex $\{a,b\}$ with a non-adjacent pair of ordered
    vertices $ab,ba$ with equal neighborhoods, and adding isolated vertices for the remaining ordered pairs.
    Since $G$ has no isolated vertices, its homomorphisms avoid these extra vertices.
    Therefore a homomorphism $G\to\widetilde
    L(H)$ is a homomorphism $G\to L(H)$ together with an independent choice of target
    ordering for each vertex of $G$, and the result follows.
\end{proof}

\begin{lem}[Edge budget]\label{lem:budget}
    For every graph $G$ and every assignment $g\colon E(G)\to\LL$,
    \[
        \left|E\left(\gex_g(G)\right)\right|\le|V(G)|.
    \]
\end{lem}
\begin{proof}
    Each gadget $\boldsymbol F=(F,(a_1,a_2),(b_1,b_2))\in\LL$ has edge
    set $\{a_1a_2,\,b_1b_2\}$, with both edges lying in an interface. Thus every edge of
    $\gex_g(G)$ is the image of an interface edge, which has an associated vertex in $G$. Since all interface edges with the same associated vertex are identified, $\gex_g(G)$ has at most one edge per vertex of $G$.
\end{proof}

\subsection{Upper bound for \texorpdfstring{$k\geq2$}{k >= 2}}\label{sec:upper}
For arbitrary gadget family $\F$ the treewidth of $\gex_g(G)$ is unbounded (a gadget may be
an arbitrarily large clique, for example). When the number of vertices of $F\in\F$ is bounded,
however, we can bound $\tw\left(\gex_g(G)\right)$. The bound is essentially that of
the \emph{blow-up} $G\boxtimes K_\mu$, obtained by replacing each vertex of $G$ by a
clique $K_\mu$ and joining adjacent cliques completely. It is easy to see, by replacing each vertex in the bags of a tree decomposition of $G$ by its $\mu$ copies, that
\[
    \tw\left(G\boxtimes K_\mu\right)\le\mu\,\tw(G)+(\mu-1).
\]
See \cite[Prop.~2.4]{DHP24} for the substitution step in the algebra of sourced graphs with distinct sourced vertices. Our gadgets have repeated interface vertices, so we provide a direct proof.

\begin{lem}[Treewidth growth]\label{lem:tw-growth}
    Let $\F$ be a gadget family with interface size $\mu$ and gadgets on at most
    $\nu$ vertices. Then for every graph $G$
    and every assignment $g\colon E(G)\to\F$,
    \[
        \tw\left(\gex_g(G)\right)\;\le\;
        \max\left\{\,\mu\,\tw(G)+(\mu-1),\;\nu-1\,\right\}.
    \]
\end{lem}
\begin{proof}
    Let $(T,(B_t)_t)$ be a tree decomposition of $G$ of width $\tw(G)$, and write $\pi$ for the quotient map defining $\gex_g(G)$.
    Replace each bag $B_t$ by $B'_t=\{\pi(u_i):u\in B_t,\ i\in[\mu]\}$.
    See \cref{fig:tw} for an illustration with $G=K_{1,3}$ and the $K_4$ gadget with two disjoint interfaces of size $2$.
    For each edge $e=uv$, let $F_e$ be its gadget copy and attach a leaf bag $\pi(V(F_e))$ to a node $t$ with $u,v\in B_t$.
    These bags cover all vertices and edges of $\gex_g(G)$.

    For a quotient vertex represented on an interface, the vertices of $G$ whose interfaces contain it form a connected subgraph: identifications occur within an interface or across an edge of $G$.
    The corresponding subtrees of $T$ therefore have connected union, and every leaf bag containing this vertex is attached to a bag containing it.
    Each internal gadget vertex appears only in its leaf bag.
    Thus the bags form a tree decomposition, with $|B'_t|\le\mu(\tw(G)+1)$ and $|\pi(V(F_e))|\le \nu$, giving
    \[
        \tw\left(\gex_g(G)\right)\le\max\left\{\,\mu\,\tw(G)+(\mu-1),\;\nu-1\,\right\}.
        \qedhere
    \]
\end{proof}

\definecolor{sfu}{HTML}{1f77b4}
\definecolor{sfv}{HTML}{2ca02c}
\definecolor{sfw}{HTML}{9467bd}

\tikzset{
  sfgv/.style={circle, fill=black, inner sep=1.8pt},
  sfprt/.style={circle, fill=black, draw=black, inner sep=0pt, minimum size=5pt},
  sfcl/.style={line width=0.7pt, black!75},
  sfspine/.style={line width=1.5pt},
  sfgedge/.style={line width=1.0pt},
  sfttl/.style={font=\small\bfseries, align=center},
}

\newcommand{\sfbagfour}[5]{\begin{scope}[on background layer]
    \fill[#1!20] (#2.center)--(#3.center)--(#4.center)--(#5.center)--cycle;
  \end{scope}
}
\newcommand{\sfbagfourline}[5]{\begin{scope}[on background layer]
    \draw[#1!75, line width=0.9pt]
      (#2.center)--(#3.center)--(#4.center)--(#5.center)--cycle;
  \end{scope}
}
\newcommand{\sfbagtwo}[3]{\begin{scope}[on background layer]
    \draw[#1!20, line width=14pt, line cap=round] (#2.center)--(#3.center);
    \draw[#1!75, line width=0.9pt, line cap=round] (#2.center)--(#3.center);
  \end{scope}
}

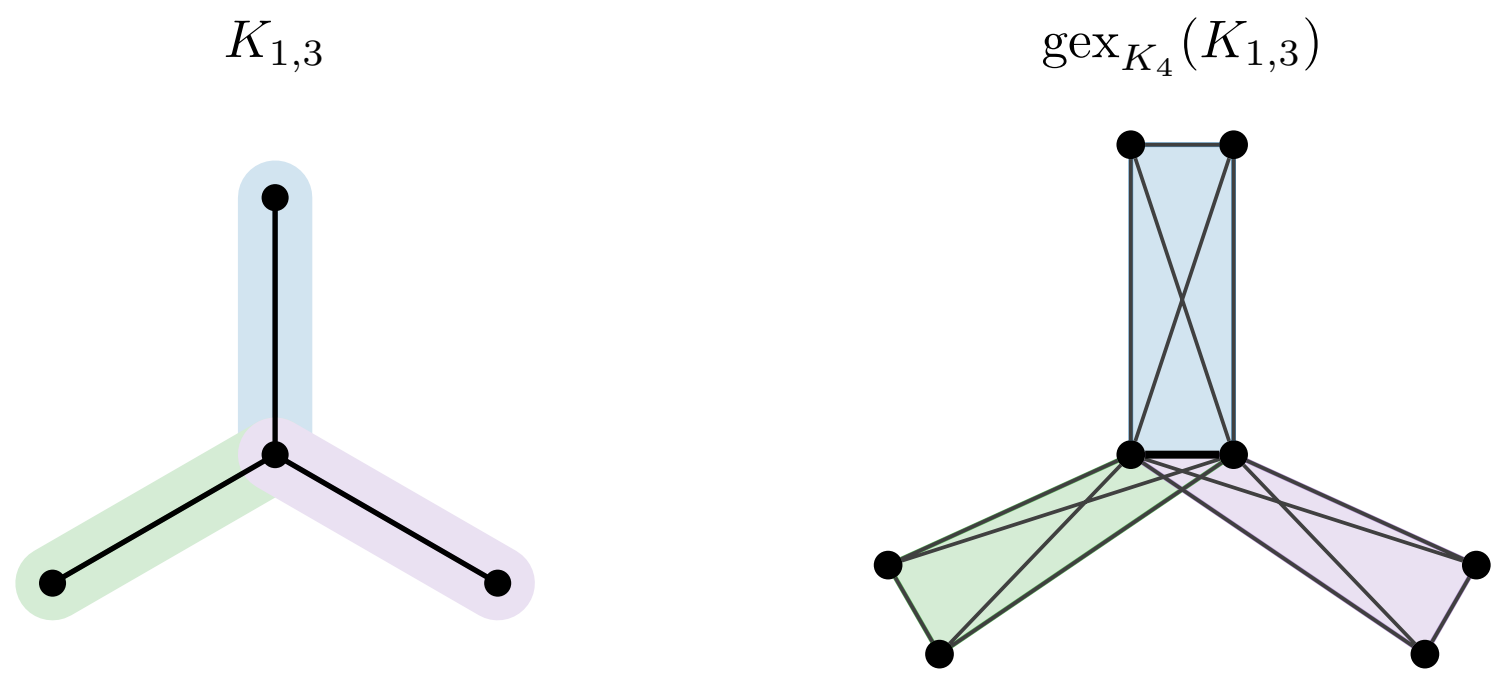
\begin{figure}
\centering
\begin{tikzpicture}

\begin{scope}[shift={(0,0)}]
  \node[sfttl] at (0,2.7) {$K_{1,3}$};
  \coordinate (c) at (0,0);
  \coordinate (u) at (90:1.7);
  \coordinate (v) at (210:1.7);
  \coordinate (w) at (330:1.7);
  \sfbagtwo{sfu}{c}{u}
  \sfbagtwo{sfv}{c}{v}
  \sfbagtwo{sfw}{c}{w}
  \draw[sfgedge] (c)--(u) (c)--(v) (c)--(w);
  \node[sfgv] at (c) {}; \node[sfgv] at (u) {};
  \node[sfgv] at (v) {}; \node[sfgv] at (w) {};
\end{scope}

\begin{scope}[shift={(6.0,0)}]
  \node[sfttl] at (0,2.7) {$\gex_{K_4}(K_{1,3})$};

  \node[sfprt] (c1) at (-0.34,0) {};
  \node[sfprt] (c2) at ( 0.34,0) {};
  \node[sfprt] (u1) at (-0.34,2.05) {};
  \node[sfprt] (u2) at ( 0.34,2.05) {};
  \node[sfprt] (v1) at ($(210:2.05)+(120:0.34)$) {};
  \node[sfprt] (v2) at ($(210:2.05)+(120:-0.34)$) {};
  \node[sfprt] (w1) at ($(330:2.05)+(60:0.34)$) {};
  \node[sfprt] (w2) at ($(330:2.05)+(60:-0.34)$) {};

  \sfbagfour{sfu}{c1}{c2}{u2}{u1}
  \sfbagfour{sfv}{c2}{c1}{v1}{v2}
  \sfbagfour{sfw}{c1}{c2}{w1}{w2}
  \sfbagfourline{sfu}{c1}{c2}{u2}{u1}
  \sfbagfourline{sfv}{c2}{c1}{v1}{v2}
  \sfbagfourline{sfw}{c1}{c2}{w1}{w2}

  \draw[sfcl] (c1)--(u1) (c1)--(u2) (c2)--(u1) (c2)--(u2) (u1)--(u2);
  \draw[sfcl] (c1)--(v1) (c1)--(v2) (c2)--(v1) (c2)--(v2) (v1)--(v2);
  \draw[sfcl] (c1)--(w1) (c1)--(w2) (c2)--(w1) (c2)--(w2) (w1)--(w2);
  \draw[sfspine] (c1)--(c2);
\end{scope}

\end{tikzpicture}
\caption{Tree decompositions of $K_{1,3}$ and $\gex_{K_{4}}(K_{1,3})=K_{1,3}\boxtimes K_2$.}
\label{fig:tw}
\end{figure}
\begin{lem}[Upper bound]\label{lem:upper-bound}
    For $k\ge1$,
    \[
        H_1\WL{2k+1}H_2\;\Longrightarrow\;L(H_1)\WL{k}L(H_2);
    \]
    this is the upper bound of \cref{thm:main-lg}.
\end{lem}
\begin{proof}
    Fix $G\in\T_k$.\footnote
    {
        Without loss of generality $G$ has no isolated vertices, as
        \cref{defn:represents} requires: writing $G=G'\sqcup\overline K_i$ with $G'$ free of
        isolated vertices, $\hom(G,X)=|V(X)|^{\,i}\hom(G',X)$, and
        $|V(L(H))|=\hom(K_2,H)/2$ is $\WLp{1}$-invariant.
    }
    For every $g\colon E(G)\to\LL$, \cref{lem:tw-growth} with $\mu=2$ (gadgets of $\LL$ have at most $3$ vertices) gives
    \[
        \tw\left(\gex_g(G)\right)\le 2\tw(G)+1\le 2k+1,
    \]
    i.e. $\Gamma_G\subseteq\T_{2k+1}$, where $\Gamma_G$ is the set of isomorphism classes
    of loopless expansions $\gex_g(G)$ over all $g$. By \eqref{eqn:decomp-iso} of \cref{lem:gadget-decomp},
    $\hom(G,L(H))$ is a linear combination of the counts $\hom(\gamma,H)$ for
    $\gamma\in\Gamma_G\subseteq\T_{2k+1}$, and $H_1\WL{2k+1}H_2$ (\cref{defn:wl}) says these counts agree. Thus $\hom(G,L(H_1))=\hom(G,L(H_2))$ for all $G\in\T_k$, i.e.
    $L(H_1)\WL{k}L(H_2)$, and the upper bound of \cref{thm:main-lg} is proved.
\end{proof}

\subsection{Lower bound for \texorpdfstring{$k\geq2$}{k >= 2}}\label{sec:lower}
We first construct a witness graph $G^*$ of treewidth $k$ with a canonical assignment expanding to $K_{2k+2}$ (\cref{lem:realization-lg,lem:witness}). By counting edges, we show that every assignment expanding $G^*$ to $K_{2k+2}$ uses only single gadgets, so $K_{2k+2}$ appears in the expansion of $\hom(G^*,L(\cdot))$ with positive coefficient (\cref{lem:positivity-lg}). Since $\tw(K_{2k+2})=2k+1$, \cref{thm:hdc} then yields $H_1\WL{2k}H_2$ with $L(H_1)\not\WL{k}L(H_2)$ (\cref{lem:separating-pair}).
The reader is encouraged to reference \cref{fig:blowup-quotient} throughout this section.
Recall the Krausz cliques $C_v=\{\,e\in E(H):v\in e\,\}$ of
\cref{defn:krausz-clique}. Throughout, $H$ and $G$ have no isolated vertices, $g\colon
E(G)\to\LL$ is a gadget assignment, and $G$ is a spanning subgraph of $L(H)$, so
$V(G)=E(H)$ and each vertex of $G$ is an edge $e=ab$ of $H$.

\definecolor{c1}{HTML}{1f77b4}
\definecolor{c2}{HTML}{ff7f0e}
\definecolor{c3}{HTML}{2ca02c}
\definecolor{c4}{HTML}{9467bd}

\tikzset{
  bqvtx/.style={circle, fill=black, inner sep=1.6pt},
  bqport/.style={circle, draw=black, line width=0.4pt, inner sep=0pt, minimum size=5.5pt},
  bqp1/.style={bqport, fill=c1}, bqp2/.style={bqport, fill=c2},
  bqp3/.style={bqport, fill=c3}, bqp4/.style={bqport, fill=c4},
  bqKvtx/.style={circle, inner sep=0pt, minimum size=7pt},
  bqK1/.style={bqKvtx, fill=c1}, bqK2/.style={bqKvtx, fill=c2},
  bqK3/.style={bqKvtx, fill=c3}, bqK4/.style={bqKvtx, fill=c4},
  bqhub/.style={draw=black, line width=0.9pt, rounded corners=2pt, inner sep=3pt},
  bqgedge/.style={line width=0.8pt},
  bqgmiss/.style={line width=0.5pt, gray!25},
  bqhedge/.style={line width=0.9pt},
  bqrel/.style={line width=0.9pt, dashed},
  bqrelmute/.style={line width=0.7pt, dashed, gray!30},
  bqslotcap/.style={rounded corners=3.5pt, draw=gray!55, line width=0.5pt, fill=gray!8},
  bqslab/.style={font=\small},
  bqtitle/.style={font=\small\bfseries, align=center},
}

\newcommand{\bqslot}[7]{\node[#4] (#1t) at (#2,{#3+0.26}) {};
  \node[#5] (#1b) at (#2,{#3-0.26}) {};
  \draw[bqhedge] (#1t)--(#1b);
  \begin{scope}[on background layer]
    \node[bqslotcap, fit=(#1t)(#1b), inner xsep=1.4pt, inner ysep=2pt] (#1cap) {};
  \end{scope}
  \node[bqslab, #7=2.5pt of #1cap] {#6};
}

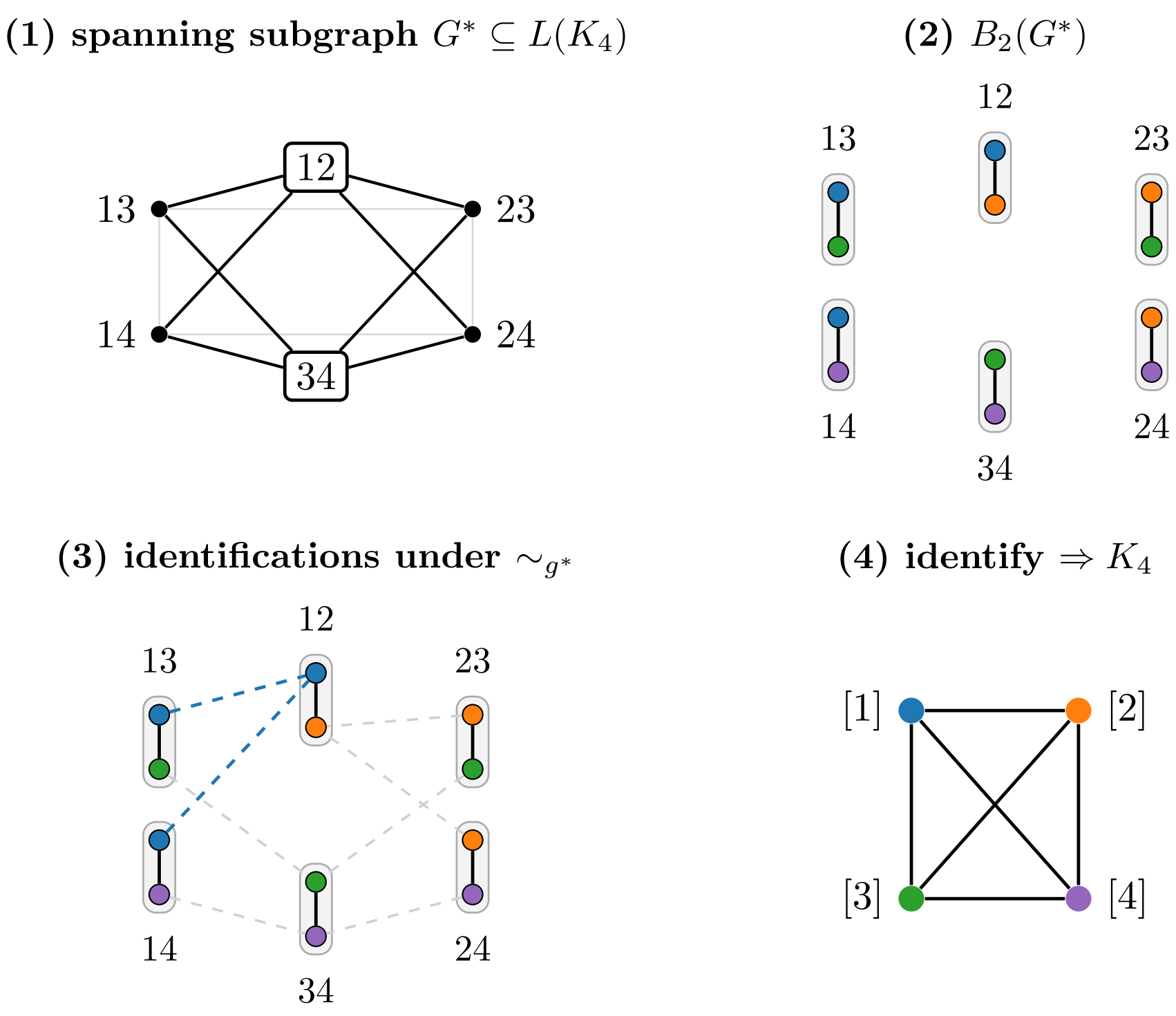
\begin{figure}
\centering
\begin{tikzpicture}

\begin{scope}[shift={(0,0)}]
  \node[bqtitle] at (1.5,2.95) {(1) spanning subgraph $G^\ast\subseteq L(K_4)$};
  \node[bqhub] (h12) at (1.5,1.7) {$12$};
  \node[bqhub] (h34) at (1.5,-0.3) {$34$};
  \node[bqvtx,label=left:$13$]  (l13) at (0,1.3) {};
  \node[bqvtx,label=left:$14$]  (l14) at (0,0.1) {};
  \node[bqvtx,label=right:$23$] (l23) at (3,1.3) {};
  \node[bqvtx,label=right:$24$] (l24) at (3,0.1) {};
  \draw[bqgmiss] (l13)--(l14);
  \draw[bqgmiss] (l23)--(l24);
  \draw[bqgmiss] (l13)--(l23);
  \draw[bqgmiss] (l14)--(l24);
  \foreach \l in {l13,l14,l23,l24}{
    \draw[bqgedge] (\l)--(h12);
    \draw[bqgedge] (\l)--(h34);
  }
\end{scope}

\begin{scope}[shift={(6.5,0)}]
  \node[bqtitle] at (1.5,2.95) {(2) $B_2(G^\ast)$};
  \bqslot{h12}{1.5}{1.6}{bqp1}{bqp2}{$12$}{above}
  \bqslot{l13}{0}{1.2}{bqp1}{bqp3}{$13$}{above}
  \bqslot{l23}{3}{1.2}{bqp2}{bqp3}{$23$}{above}
  \bqslot{h34}{1.5}{-0.4}{bqp3}{bqp4}{$34$}{below}
  \bqslot{l14}{0}{0.0}{bqp1}{bqp4}{$14$}{below}
  \bqslot{l24}{3}{0.0}{bqp2}{bqp4}{$24$}{below}
\end{scope}

\begin{scope}[shift={(0,-5)}]
  \node[bqtitle] at (1.5,2.95) {(3) identifications under $\sim_{g^*}$};
  \bqslot{h12}{1.5}{1.6}{bqp1}{bqp2}{$12$}{above}
  \bqslot{l13}{0}{1.2}{bqp1}{bqp3}{$13$}{above}
  \bqslot{l23}{3}{1.2}{bqp2}{bqp3}{$23$}{above}
  \bqslot{h34}{1.5}{-0.4}{bqp3}{bqp4}{$34$}{below}
  \bqslot{l14}{0}{0.0}{bqp1}{bqp4}{$14$}{below}
  \bqslot{l24}{3}{0.0}{bqp2}{bqp4}{$24$}{below}
  \draw[bqrelmute] (l23t)--(h12b);
  \draw[bqrelmute] (l24t)--(h12b);
  \draw[bqrelmute] (l13b)--(h34t);
  \draw[bqrelmute] (l23b)--(h34t);
  \draw[bqrelmute] (l14b)--(h34b);
  \draw[bqrelmute] (l24b)--(h34b);
  \draw[bqrel, c1] (l13t)--(h12t);
  \draw[bqrel, c1] (l14t)--(h12t);
\end{scope}

\begin{scope}[shift={(6.5,-5)}]
  \node[bqtitle] at (1.5,2.95) {(4) identify $\Rightarrow K_4$};
  \node[bqK1, label=left:{$[1]$}]  (v1) at (0.7,1.5) {};
  \node[bqK2, label=right:{$[2]$}] (v2) at (2.3,1.5) {};
  \node[bqK3, label=left:{$[3]$}]  (v3) at (0.7,-0.3) {};
  \node[bqK4, label=right:{$[4]$}] (v4) at (2.3,-0.3) {};
  \foreach \a/\b in {v1/v2,v1/v3,v1/v4,v2/v3,v2/v4,v3/v4}{
    \draw[bqhedge] (\a)--(\b);
  }
\end{scope}

\end{tikzpicture}
\caption{Example illustrating the canonical map $\gex_{g^*}(\cdot)$ sending the spanning subgraph $G^*\subseteq L(K_4)$ to $K_4$. $G^*$ is the witness graph for $K_4$ produced by the matching $\{12,34\}$. The color of a vertex records its image under the label map $\pi$. Note that $K_4$ lies outside the setting of the main theorem; this figure is for illustration only.}
\label{fig:blowup-quotient}
\end{figure}

\begin{defn}[$B_2$, label map]\label{defn:b2}
    Let $B_2$ denote the operation sending each vertex of a graph to a labeled copy of
    $K_2$. In particular, for $G$ as above,
    \[
        B_2(G)=\bigsqcup_{e\in V(G)}B_2(e)
    \]
    where the two vertices of $B_2(e)\cong K_2$ are labeled by the endpoints of $e$.
    Write
    \[
        \pi\colon V(B_2(G))\to V(H)
    \]
    for the \emph{label map}, sending each vertex of $B_2(G)$ to its label.
\end{defn}

\noi Since the edges of gadgets in $\LL$ are contained in their interfaces, the vertices
and edges of $\gex_g(G)$ are exactly those of $B_2(G)$ modulo the identifications of the
$\gex_g$ operation. Thus, denoting the equivalence relation generated by those
identifications as $\sim_g$, we have
\[
    \gex_g(G)=B_2(G)/\!\sim_g
\]
(with parallel edges collapsed). Write $\pi_g\colon V(B_2(G))\to V(\gex_g(G))$ for the
associated quotient map.

\begin{defn}[Canonical assignment]\label{defn:canonical-lg}
    Let $ef\in E(G)$. Since $G$ is a subgraph of $L(H)$ we have $ef\in E(L(H))$, and thus
    there is a unique vertex $v\in V(H)$ with $ef$ in the Krausz clique $C_v$. The
     \emph{canonical assignment} $g^*$ assigns to $ef$ the gadget identifying the
    unique vertices of $B_2(e)$ and $B_2(f)$ with label $v$.
\end{defn}

\begin{lem}[Realization]\label{lem:realization-lg}
    With $G$ a spanning subgraph of $L(H)$ and $g^*$ the canonical assignment  above,
    \[
        \gex_{g^*}(G)\cong H
        \quad\Longleftrightarrow\quad
        G[C_v]\text{ is connected for every }v\in V(H).
    \]
\end{lem}

\begin{proof}
    Since $g^*$ only identifies vertices of $B_2(G)$ with the same label, $\pi$ is
    constant on classes of $\sim_{g^*}$ and factors as
    \[
        \pi=\bar\pi_{g^*}\circ\pi_{g^*}
    \]
    with $\bar\pi_{g^*}\colon V(\gex_{g^*}(G))\to V(H)$ sending each class of
    $\sim_{g^*}$ to the common label of its members. If $\pi(u)\neq\pi(u')$ then
    $[u]\neq[u']$, so $\bar\pi_{g^*}$ separates classes of distinct label and
    $\gex_{g^*}(G)$ is loopless.

    Fix $v\in V(H)$. The label-$v$ vertex of $B_2(e)$ and that of $B_2(f)$ are identified
    by $g^*$ exactly when $ef\in E(G[C_v])$. Passing to $\sim_{g^*}$, two vertices
    of label $v$ are equivalent iff their edges lie in one connected component of
    $G[C_v]$, i.e.
    \[
        \bar\pi_{g^*}^{-1}(v)=\{\text{connected components of }G[C_v]\}.
    \]
    Note $H$ has no isolated vertices so $C_v\neq\varnothing$ and $\bar\pi_{g^*}$ is onto;
    the label classes partition $V(B_2(G))$, so the fibers partition
    $V(\gex_{g^*}(G))$. Finally, each vertex $e=ab\in V(G)$ contributes an edge of
    $B_2(G)$ whose image in $\gex_{g^*}(G)$ is $[a][b]$; as $G$ spans $L(H)$,
    $\bar\pi_{g^*}$ induces a bijection $E(\gex_{g^*}(G))\to E(H)$.

    To complete the proof, in one direction, if every $G[C_v]$ is connected, each fiber $\bar\pi_{g^*}^{-1}(v)$
    is a single class, so $\bar\pi_{g^*}$ is a bijection on vertices. Since it is also a bijection on edges it is an isomorphism, and $\gex_{g^*}(G)\cong H$.

    In the other direction, if some $G[C_v]$ is disconnected, then $\bar\pi_{g^*}^{-1}(v)$ has
    at least two classes, so $|V(\gex_{g^*}(G))|>|V(H)|$ and
    $\gex_{g^*}(G)\not\cong H$.
\end{proof}

\begin{defn}[Witness graph]\label{defn:witness-lg}
    Let $n=2k+2$ and fix a perfect matching $M=\{m_1,\dots,m_{k+1}\}$ of $K_n$. For
    $v\in V(K_n)$ let $m(v)\in M$ be the matching edge incident to $v$. The
    \emph{witness graph} $G^*$ is the spanning subgraph of $L(K_n)$ whose edges join
    each matching edge to every other edge of its Krausz clique. That is,
    \[
        V(G^*)=E(K_n),\qquad
        E(G^*)=\left\{\,m(v)\,e \;:\; v\in V(K_n),\ e\in C_v,\ e\neq m(v)\,\right\}.
    \]
    Each such pair $m(v),e$ shares the vertex $v$, so $G^*\subseteq L(K_n)$.
\end{defn}

\begin{lem}[Witness graph]\label{lem:witness}
    For $k\ge2$, let $g^*$ be the canonical assignment of $G^*\subseteq L(K_{2k+2})$.
    Then
    $\gex_{g^*}(G^*)\cong K_{2k+2}$ and $\tw(G^*)=k$, so
    \[
        \tw(\gex_{g^*}(G^*))=2\tw(G^*)+1.
    \]
\end{lem}

\begin{proof}
    For each $v\in V(K_{2k+2})$, the Krausz clique $C_v$ induces a star in $G^*$ with
    center $m(v)$, which is connected. By \cref{lem:realization-lg},
    $\gex_{g^*}(G^*)\cong K_{2k+2}$, and thus has treewidth $2k+1$.

    For $\tw(G^*)\le k$, we exhibit a tree decomposition. Create a central bag
    $B_0=M=\{m_1,\dots,m_{k+1}\}$ of the $k+1$ matching edges, and for each non-matching
    edge $e=uv$ a bag $B_e=\{m(u),m(v),e\}$ adjacent to $B_0$. Every vertex appears in
    some bag; the edges $m(u)e$ and $m(v)e$ of $G^*$ both lie in $B_e$; each matching
    edge $m_i$ occupies $B_0$ together with all bags $B_e$ containing it, a connected
    subtree; and each non-matching edge occupies a single bag. The width is
    $\max(k+1,3)-1=k$ for $k\ge2$, so $\tw(G^*)\le k$.

    Finally, \cref{lem:tw-growth} gives
    \[
        2k+1=\tw(\gex_{g^*}(G^*))\le 2\,\tw(G^*)+1,
    \]
    so $\tw(G^*)\ge k$ (the second term of the maximum in \cref{lem:tw-growth} is at most $2<2k+1$). Thus $\tw(G^*)=k$ and the bound of \cref{lem:tw-growth} is
    achieved.
\end{proof}

\begin{rem}\label{rem:k1-tree}
    At $k=1$, the construction in \cref{defn:witness-lg} gives $G^*\cong K_{2,4}$, of treewidth $2$, with canonical expansion $K_4$.
    Its treewidth is therefore too large to witness a failure of transfer at $k=1$, consistent with $\tau_L(1)=1$.
\end{rem}

\begin{lem}[Positivity]\label{lem:positivity-lg}
    $\beta(G^*,K_{n})\geq 1$.
\end{lem}

\begin{proof}
    By \cref{lem:budget}, $|E(\gex_g(G^*))|\le|V(G^*)|$ for every assignment $g$,
    and since $G^*$ spans $L(K_n)$ we have $|V(G^*)|=|E(K_n)|=\binom n2$.

    Now suppose $\gex_g(G^*)\cong K_n$. Then it has exactly $\binom n2$ edges, so the
    budget is tight: no edges of $B_2(G^*)$ are identified by $g$, and thus $g$
    assigns no double gadget $\boldsymbol F_\parallel$ or $\boldsymbol F_\times$. Hence
    \[
        \prod_{e\in E(G^*)}\alpha(g(e))=(+1)^{|E(G^*)|}=+1,
    \]
    and therefore
    \[
        \beta(G^*,K_n)
        =\sum_{g:\,\gex_g(G^*)\cong K_n}\ \prod_{e\in E(G^*)}\alpha(g(e))
        =\#\{g:\gex_g(G^*)\cong K_n\}\ge1,
    \]
    with the canonical $g^*$ being one such assignment (\cref{lem:witness}).
\end{proof}

\begin{lem}[Separating pair]\label{lem:separating-pair}
    For every $k\ge2$ there exist $H_1\WL{2k}H_2$ with $L(H_1)\not\WL{k}L(H_2)$; this is the lower bound of \cref{thm:main-lg}.
\end{lem}

\begin{proof}
    By \eqref{eqn:decomp-iso} applied to the witness $G^*$,
    \[
        \hom\left(G^*,L(H)\right)
        =2^{-|V(G^*)|}\sum_{\gamma\in\Gamma}\beta(G^*,\gamma)\cdot\hom(\gamma,H),
    \]
    with $\Gamma=\{\gamma\in\G:\beta(G^*,\gamma)\neq0\}$, a finite set of pairwise non-isomorphic graphs. By \cref{lem:positivity-lg} $K_n\in \Gamma$, and $\tw(K_n)=2k+1>2k$, so
    $\Gamma\not\subseteq\T_{2k}$.
    Now suppose for contradiction that $\hom(G^*,L(\cdot))$ is
    $\WL{2k}$-invariant, i.e.
    \[
        H_1\WL{2k}H_2\;\Longrightarrow\;\hom(G^*,L(H_1))=\hom(G^*,L(H_2)).
    \]
    Clearing the
    constant $2^{-|V(G^*)|}$, this is exactly the hypothesis of \cref{thm:hdc} for
    $\left(\Gamma,\beta(G^*,\cdot)\right)$ and the class $\T_{2k}$, whose conclusion is $\Gamma\subseteq\T_{2k}$. This contradicts $K_n\in \Gamma$, and hence there exist $H_1\WL{2k}H_2$ with
    \[
        \hom\left(G^*,L(H_1)\right)\neq\hom\left(G^*,L(H_2)\right).
    \]
    Then $\tw(G^*)=k$ implies $L(H_1)\not\WL{k}L(H_2)$, proving this lemma and the lower bound of \cref{thm:main-lg}.
\end{proof}

\subsection{Special case \texorpdfstring{$k=1$}{k = 1}}\label{sec:k1}

For connected graphs the transfer result for the line graph at $k=1$ also follows from covering theory, as \cref{rem:covering} sketches. We record a proof supplied by our framework for more general gadget families that we will also need in \cref{sec:rs}. It relies on a treewidth
bound for expansions of forests that is sharper than \cref{lem:tw-growth}.

We let $G$ be a tree with at least one edge and $g\colon E(G)\to\LL$ be an assignment.
Forests follow componentwise, since an isolated vertex expands to two isolated vertices.
Recall from \cref{defn:gex} that
$\gex_g(G)$ is composed of the interface vertices and edges
\[
    \left(\bigsqcup_{x\in V(G)}\{x_1,x_2\}\right)/\sim_g,\qquad\left(\bigsqcup_{x\in V(G)}x_1x_2\right)/\sim_g.
\]
Each tree edge $e=xy$ contributes, via the gadget
$g(e)$, one or two \emph{identification edges} recording the identifications $x_i=y_j$
that $g(e)$ prescribes: a single gadget $\boldsymbol F_{ij}$ contributes the one edge
$x_iy_j$; a double gadget contributes the two edges $\{x_1y_1,x_2y_2\}$ or
$\{x_1y_2,x_2y_1\}$. Let $R$ be the \emph{identification graph} on $V(G)\times[2]$
carrying these edges, so that the components of $R$ are the classes of $\sim_g$, i.e. the
vertices of $\gex_g(G)$.

\begin{lem}\label{lem:R-forest}
    $R$ is a forest, and $x_1,x_2$ lie in distinct components of $R$ for every $x\in V(G)$.
\end{lem}

\begin{proof}
    Suppose for contradiction that $W\subseteq R$ is either a simple cycle or a simple
    path from $x_1$ to $x_2$ for some $x\in V(G)$. Projecting each interface vertex $x_i$
    to $x\in V(G)$ sends $W$ to a closed walk in $G$ (both endpoints of a path
    $x_1\!-\!x_2$ project to $x$). Since $G$ is a tree this walk contains an immediate
    backtrack: two consecutive edges of $W$ projecting to the same tree edge of $G$ and thus
    belonging to the same gadget. This contradicts the fact that, by inspection, the identification graph $R$ is simple (no loops or parallel edges), and the identification edges of any one gadget are vertex-disjoint.
\end{proof}

\begin{prop}[Forest expansion]\label{prop:forest-comp}
    For every forest $G$ and every assignment $g\colon E(G)\to\LL$, $\gex_g(G)$ is a
    forest.
\end{prop}

\begin{proof}
    Write $\pi_g$ for the quotient map onto $V(\gex_g(G))$. Take $G$ itself as the
    decomposition tree, with bag $B_x=\{\pi_g(x_1),\pi_g(x_2)\}$ at $x\in V(G)$; by
    \cref{lem:R-forest} these are two distinct vertices. Every vertex of $\gex_g(G)$ is the
    class of some interface vertex and so lies in some bag. Every edge of $\gex_g(G)$ is the
    image of some interface edge $x_1x_2$, since the gadgets of $\LL$ carry no other edges,
    and so lies in $B_x$. Finally, for a fixed class $c$, the set $\{x: c\in B_x\}$ is the
    projection to $V(G)$ of a component of $R$, which is connected in $G$ because $R$
    carries edges only along edges of $G$; a connected subgraph of a tree is a subtree. Thus
    $(G,(B_x)_x)$ is a tree decomposition of width $1$, and $\gex_g(G)$ is a forest.
\end{proof}

\begin{cor}\label{cor:k1-transfer}
    For all graphs $H_1,H_2$, $H_1\WL{1}H_2$ implies $L(H_1)\WL{1}L(H_2)$; this is
    the case $k=1$ of \cref{thm:main-lg}.
\end{cor}
\begin{proof}
    Write a forest $G$ as $G'\sqcup\overline K_i$, where $G'$ has no isolated vertices.
    Then $\hom(G,L(H))=(\hom(K_2,H)/2)^i\hom(G',L(H))$.
    By \eqref{eqn:decomp-iso}, the last factor is a rational combination of the counts $\hom(\gamma,H)$ with each
    $\gamma=\gex_g(G')$ a forest by \cref{prop:forest-comp}, i.e. $\gamma\in\T_1$.
    If $H_1\WL{1}H_2$ these counts and $\hom(K_2,H)$ agree, so $\hom(G,L(H_1))=\hom(G,L(H_2))$ for every forest
    $G$, i.e. $L(H_1)\WL{1}L(H_2)$, and the case $k=1$ of \cref{thm:main-lg} is proved.
\end{proof}

\begin{rem}\label{rem:covering}
    For connected graphs, the corollary also follows from covering theory. Connected graphs $H_1,H_2$ satisfy $H_1\WL{1}H_2$ iff they have equal order and a common finite cover~\cite{KV15}; by Leighton's theorem~\cite{Leighton82}, the latter holds iff they have isomorphic universal covers.
    So let $H_1\WL{1}H_2$ be connected, and let $X\to H_1$ and $X\to H_2$ be finite covers. Since $H_1$ and $H_2$ have equal order, both covers are $d$-sheeted. Then the covers $X\to H_1$ and $X\to H_2$ induce $d$-sheeted covers $L(X)\to L(H_1)$ and $L(X)\to L(H_2)$. Thus $L(X)$ is a common finite cover of $L(H_1)$ and $L(H_2)$, which are connected and have the same order $|E(X)|/d$, and so $L(H_1)\WL{1}L(H_2)$.
    Iterating yields $L^{j}(H_1)\WL{1}L^{j}(H_2)$ for all $j\geq1$.
\end{rem}

\section{\texorpdfstring{$(r,S)$}{(r,S)}-line graphs: Proof of \texorpdfstring{\cref{thm:main-rs}}{the (r,S)-line-graph transfer theorem}}\label{sec:rs}

Letting $H$ be a graph on $[n]$, we represent each $r$-clique by a tuple of its
vertices $P=(p_1,\hdots,p_r)\in[n]^r$ and denote its underlying set by $\{P\}$. Recall from \cref{defn:rs} that for $r\ge2$ and $S\subseteq\{0,1,\dots,r-1\}$ the
$(r,S)$-line graph $\LrS(H)$ has the $r$-cliques of $H$ as vertices, with $P\sim Q$ whenever
$|\{P\}\cap \{Q\}|\in S$.

Fix $r\ge2$ and $S\subseteq\{0,1,\dots,r-1\}$ for the rest of this section. Write $\Lr=\LrS$ for the $(r,S)$-line graph operation and $\LL_r=\LL_{r,S}$ for its gadget family. The upper bounds of \cref{cor:tightness,cor:k1-rs} hold for every such $S$. The lower bounds for $k\ge2$ require the additional hypothesis $\{1\}\subseteq S\subseteq\{1,\dots,r-1\}$: $1\in S$ is required for the witness construction of \cref{prop:witness-realizes}, and $0\not\in S$ for the non-vanishing coefficients of \cref{lem:witness-coeff,prop:cycle-witness}. As in \cref{subsec:gadgets-lg}, we go through ordered $r$-cliques to expand $\hom(G,\Lr(H))$.

\begin{defn}[Ordered $(r,S)$-line graph]\label{defn:ordered-rs}
    For $P\in[n]^r$ let $\kappa_H(P)=\prod_{1\le i<j\le r}H_{p_i,p_j}$ be the \emph{ordered
    clique indicator}, equal to $1$ when $\{P\}$ is an $r$-clique of $H$ and $0$ otherwise. With
    \[
        A_S(P,Q)=\bone\!\left[\,|\{P\}\cap\{Q\}|\in S\,\right],
    \]
    the \emph{ordered $(r,S)$-line graph} $\widetilde\Lambda_r(H)$ is the undirected
    graph on $[n]^r$ with adjacency
    \[
        \widetilde\Lambda_r(H)_{P,Q}=\kappa_H(P)\,\kappa_H(Q)\,A_S(P,Q).
    \]
\end{defn}

\subsection{The \texorpdfstring{$(r,S)$}{(r,S)} gadget family}\label{subsec:gadgets-rs}

\noindent Fix a graph $H$ on $[n]$ and $P,Q\in V(\widetilde{\Lambda}_r(H))=[n]^r$ with distinct entries throughout this section, and let $M$ denote a (partial) matching $M\subseteq[r]\times[r]$; in particular we record the data of $\{P\}\cap\{Q\}$ by the \emph{coincidence matching} $M^*\subseteq[r]\times[r]$ which has size $|\{P\}\cap\{Q\}|$ since $P,Q$ have no repeated entries. We suppress the dependence on $P,Q$ and define the indicators
\[
    E_M=\prod_{(i,j)\in M}\bone[p_i=q_j]=\bone[M\subseteq M^*],\qquad \delta_M=\bone[M=M^*],
\]
and observe
\begin{equation}\label{eqn:AS-delta}
    A_S=\sum_{M\;:\;|M|\in S}\delta_M,
    \qquad
    E_M=\sum_{M'\supseteq M}\delta_{M'}.
\end{equation}

\begin{lem}[Coefficients]\label{lem:coefficients}
    For $P,Q\in[n]^r$ with distinct entries,
    \begin{equation}\label{eqn:mobius}
        A_S(P,Q)=\sum_{M}\alpha_{|M|}\;E_M(P,Q),\qquad \alpha_t=\sum_{\substack{s\leq t \\ s\in S}}(-1)^{t-s}\binom ts,
    \end{equation}
    where $M$ ranges over all matchings $M\subseteq [r]\times [r]$.
\end{lem}

\begin{proof}
    Consider the poset of matchings ordered by inclusion. The interval $[M,M']$ is the Boolean lattice on $M'\setminus M$, with M\"obius function
    $\mu(M,M')=(-1)^{|M'|-|M|}$. Inverting the second identity of \eqref{eqn:AS-delta},
    \[
        \delta_M=\sum_{M'\supseteq M}(-1)^{|M'|-|M|}\;E_{M'}.
    \]
    Then substituting into the first identity of \eqref{eqn:AS-delta},
    \[
        A_S=\sum_{M\;:\;|M|\in S}\ \sum_{M'\supseteq M}(-1)^{|M'|-|M|}\;E_{M'}
        =\sum_{M'}\sum_{\substack{M\subseteq M'\\ |M|\in S}}(-1)^{|M'|-|M|}E_{M'}.
    \]
    Then letting $t=|M'|$, there are $\binom ts$ sub-matchings $M\subseteq M'$ of size $s$, and the coefficient identity of \eqref{eqn:mobius} follows.
\end{proof}

\noindent For $r=2$, $S=\{1\}$ this gives $(\alpha_0,\alpha_1,\alpha_2)=(0,+1,-2)$,
recovering the coefficients of $\LL$ (\cref{prop:ie-lg}).
As a consistency check, the coefficients recover the adjacency rule of $\Lambda_{r,S}$ on cliques:

\begin{cor}[Adjacency recovery]\label{cor:adjacency-recovery}
    For every $0\le j\le r$,
    \[
        \sum_{t=0}^{j}\binom jt\,\alpha_t=\bone[j\in S].
    \]
\end{cor}

\begin{proof}
    Substitute $\alpha_t=\sum_{s\in S,\,s\le t}(-1)^{t-s}\binom ts$ and exchange summation:
    \[
        \sum_{t=0}^{j}\binom jt\alpha_t
        =\sum_{\substack{s\in S\\ s\le j}}\ \sum_{t=s}^{j}\binom jt\binom ts(-1)^{t-s}.
    \]
    Reorganizing by the identity $\binom jt\binom ts=\binom js\binom{j-s}{t-s}$ and re-indexing the inner sum by $u=t-s$ then yields
    \[
        \binom js\sum_{u=0}^{j-s}\binom{j-s}{u}(-1)^{u}=\binom js(1-1)^{j-s}=\bone[j=s]. \qedhere
    \]
\end{proof}

\subsection{The case \texorpdfstring{$k\geq2$}{k >= 2}}
\label{subsec:kge2-rs}

We first build a gadget family representing $\Lr$
(\cref{defn:family-rs,lem:clique-gadget,lem:duplication-rs}). We then describe when a
spanning subgraph of $\Lr(K_n)$ expands to $K_n$ (\cref{lem:realization-rs}), exhibit a
graph of treewidth $k$ that does so (\cref{defn:witness-rs,lem:witness-tw}) with nonzero
coefficient (\cref{lem:witness-coeff}), and conclude with a lower bound for
\cref{thm:main-rs} in \cref{cor:tightness}.

\begin{defn}[The $(r,S)$-line graph gadget family]\label{defn:family-rs}
    The \emph{double clique gadget} $\boldsymbol F_0=(F_0,\boldsymbol{a},\boldsymbol{b})$ is the
    bi-labeled graph of interface size $r$ on $\{a_i\}\cup\{b_j\}$ whose edges form cliques on $\{a_i\}$ and $\{b_j\}$ for all $i,j\in[r]$. For a matching $M\subseteq[r]\times[r]$, the gadget
    $\boldsymbol F_M$ is obtained from $\boldsymbol F_0$ by identifying $(a_i{=}b_j)$
    for each $(i,j)\in M$. The \emph{$(r,S)$-line graph gadget family} is
    \[
        \LL_r=\{\,\boldsymbol F_M :  \alpha_{|M|}\neq0\,\},\qquad \alpha(\boldsymbol F_M)=\alpha_{|M|},
    \]
    of interface size $r$.
\end{defn}

\noi When $0\not\in S$, $\alpha_0=0$, so
$\boldsymbol F_0\not\in\LL_r$, and every gadget of $\LL_r$ identifies some pair of vertices. For $r=2$, $S=\{1\}$ the matchings of size $1$ give the four single gadgets
$\boldsymbol F_{11},\boldsymbol F_{12},\boldsymbol F_{21},\boldsymbol F_{22}$ ($\alpha=+1$)
and those of size $2$ the doubles $\boldsymbol F_\parallel,\boldsymbol F_\times$
($\alpha=-2$), recovering $\LL$ (\cref{defn:family-lg}).

\begin{lem}[Clique-gadget matrices]\label{lem:clique-gadget}
    For every matching $M\subseteq[r]\times[r]$ and every graph $H$ on $[n]$,
    \[
        T^{\boldsymbol F_M\to H}_{P,Q}=\kappa_H(P)\,\kappa_H(Q)\,E_M(P,Q),\qquad
        \widetilde\Lambda_r(H)=\sum_{\boldsymbol F\in\LL_r}\alpha(\boldsymbol F)\,T^{\boldsymbol F\to H},
    \]
    where $\kappa_H$ is the ordered clique indicator of \cref{defn:ordered-rs}.
\end{lem}

\begin{proof}
    Write $\boldsymbol{F}_M=(F_M,\boldsymbol{a},\boldsymbol{b})$. Since every vertex of $\boldsymbol{F}_M$ is an interface vertex, if there exists a homomorphism $\varphi:V(\boldsymbol{F}_M)\to[n]$ with $\varphi:\boldsymbol{a}\mapsto P$ and $\boldsymbol{b}\mapsto Q$, it is uniquely determined, and thus $T_{P,Q}^{\boldsymbol{F}_M\to H}\in\{0,1\}$. Such a homomorphism exists iff $P,Q$ are consistent with the identifications of $M$ and send edges of $F_M$ to edges of $H$, i.e. $\kappa_H(P)=\kappa_H(Q)=E_M(P,Q)=1$, and the first equality is shown.
    Now $\alpha_{|M|}=0$ for $\boldsymbol{F}_M\not\in\LL_r$, so
    \[
        \sum_{\boldsymbol F\in\LL_r}\alpha(\boldsymbol F)\,T^{\boldsymbol F\to H}_{P,Q}
        =\kappa_H(P)\kappa_H(Q)\sum_M\alpha_{|M|}E_M(P,Q).
    \]
    If either $P$ or $Q$ have repeated entries, then $\kappa_H(P)\kappa_H(Q)$ vanishes. Otherwise,
    \cref{lem:coefficients} gives $\sum_M\alpha_{|M|}E_M(P,Q)=A_S(P,Q)$, and the right hand side is $\kappa_H(P)\kappa_H(Q)A_S(P,Q)=\widetilde\Lambda_r(H)_{P,Q}$.
\end{proof}

\begin{lem}[Duplication factor]\label{lem:duplication-rs}
    For every graph $H$ and every graph $G$ without isolated vertices,
    \[
        \hom\!\left(G,\widetilde\Lambda_r(H)\right)=(r!)^{|V(G)|}\,\hom\!\left(G,\Lr(H)\right),
    \]
    so $\LL_r$ represents $\Lr$ with duplication factor $d=r!$.
\end{lem}

\begin{proof}
    The argument is that of \cref{lem:duplication-lg} with $r!$ orderings in place of $2$.
    Tuples that do not represent $r$-cliques are isolated in $\widetilde\Lambda_r(H)$ and cannot occur in the image of a homomorphism from $G$.
\end{proof}

Fix $n$ and take $H=K_n$ on $[n]$ (where every $r$-subset is an
$r$-clique). Let $\Omega$ be a set of $r$-subsets of $[n]$ and $G^*$ be a spanning subgraph
of $\Lr(K_n)[\Omega]$, i.e.
\[
    V(G^*)=\Omega,\qquad E(G^*)\subseteq\{\,PQ:P,Q\in\Omega,\ |P\cap Q|\in S\,\}.
\]
For $v\in[n]$ write
$C_v=\{\,P\in\Omega:v\in P\,\}$ for the \emph{Krausz class} of $v$; at $r=2$ these are the
Krausz cliques of \cref{subsec:krausz}. \cref{defn:br} and \cref{defn:canonical-rs} directly generalize \cref{defn:b2} and \cref{defn:canonical-lg}, while \cref{lem:realization-rs} has an additional condition beyond \cref{lem:realization-lg}.

\begin{defn}[$B_r$, label map]\label{defn:br}
    Let $B_r$ denote the operation sending each vertex of a graph to a labeled copy of
    $K_r$. For $G^*$ as above,
    \[
        B_r(G^*)=\bigsqcup_{P\in\Omega}B_r(P),
    \]
    where the $r$ vertices of $B_r(P)=K_r$ are labeled bijectively by the elements of
    $P$. Write $\pi\colon V(B_r(G^*))\to[n]$ for the \emph{label map}, sending each
    vertex to its label.
\end{defn}

\noindent Suppose $G^*$ has no isolated vertices.
Since every gadget $\boldsymbol{F}_M\in\LL_r$ consists of the two $r$-cliques on its
interfaces, the vertices and edges of $\gex_g(G^*)$ are those of $B_r(G^*)$ modulo
the identifications of $\gex_g$. Thus, writing $\sim_g$ for the equivalence relation generated by those relations,
\[
    \gex_g(G^*)=B_r(G^*)/\!\sim_g
\]
(with parallel edges collapsed), with quotient map $\pi_g\colon V(B_r(G^*))\to
V(\gex_g(G^*))$.

\begin{defn}[Canonical assignment]\label{defn:canonical-rs}
    The \emph{canonical assignment} $g^*$ assigns to each edge $PQ\in E(G^*)$ the
    gadget identifying the $v$-labeled vertex of $B_r(P)$ with
    the $v$-labeled vertex of $B_r(Q)$ for every $v\in P\cap Q$.
    This is an assignment into $\LL_r$ provided $\alpha_{|P\cap Q|}\neq0$ for every edge $PQ\in E(G^*)$; for the witness graph of \cref{defn:witness-rs}, every edge has $|P\cap Q|=1$ and $\alpha_1=1$.
\end{defn}

\begin{lem}[Realization]\label{lem:realization-rs}
    Let $G^*$ be a spanning subgraph of $\Lr(K_n)[\Omega]$ without isolated vertices, with canonical assignment $g^*$ as above.
    Then
    \[
        \gex_{g^*}(G^*)\cong K_n
        \quad\Longleftrightarrow\quad
        \begin{cases}
            \text{(i)\,\, $G^*[C_v]$ is connected for every $v\in[n]$,}\\[2pt]
            \text{(ii) every pair of $[n]$ lies in some $P\in\Omega$.}
        \end{cases}
    \]
\end{lem}

\begin{proof}
    The proof is the same as \cref{lem:realization-lg} aside from the additional condition $(ii)$. In the forward direction, note that two classes with the same label are never adjacent, since the vertices of each $B_r(P)$ carry distinct labels. So $\gex_{g^*}(G^*)\cong K_n$ implies each of the $n$ labels has exactly one class, which gives $(i)$. For $(ii)$, notice each $P=\{p_1,\dots,p_r\}\in\Omega$
    contributes the $\binom r2$ clique edges among $[p_1],\dots,[p_r]$, so $uv$ is an edge
    of $\gex_{g^*}(G^*)$ iff some $P\in\Omega$ contains both $u$ and $v$. At $r=2$
    this is automatic from $G$ spanning $L(H)$.
\end{proof}

From \cref{defn:witness-rs} through \cref{lem:witness-coeff}, assume $\{1\}\subseteq S\subseteq\{1,\dots,r-1\}$.
The witness is built so that its treewidth is controlled by a single large bag and every realizing assignment uses only single gadgets, giving us the treewidth and non-vanishing coefficient respectively.

\begin{defn}[Witness graph]\label{defn:witness-rs}
Let $n=rm$ with $m\ge r$ and fix a partition $\mathcal P=\{P_1,\dots,P_m\}$ of $K_n$ on $[n]$ into
$m$ \emph{partition cliques} of size $r$. Then fix a second family $\mathcal Q=\{Q_1,\dots,Q_t\}$ of \emph{cross cliques} satisfying
\begin{itemize}
    \item Every cross clique contains at most one vertex from each partition clique.
    \item For every pair $\{u,v\}\in{\binom{[n]}{2}}$ with $u,v$ in different partition cliques, $\{u,v\}$ lies in some cross clique.
\end{itemize}
Such a family exists since $m\ge r$: for instance, all $r$-subsets of $[n]$ meeting $r$ distinct partition cliques in one vertex each.

\noindent Put $\Omega=\mathcal P\cup\mathcal Q$, and let
$G^*$ be the spanning subgraph of $\Lr(K_n)[\Omega]$ joining each partition clique $P\in\mathcal P$ to all cross cliques $Q\in\mathcal Q$ it intersects, i.e.
\[
    E(G^*)=\{PQ : (P,Q)\in\mathcal P\times\mathcal Q,P\cap Q\neq\varnothing\,\}.
\]
\end{defn}

The hypothesis $1\in S$ ensures $PQ$ is an edge of $\Lr(K_n)$. At $r=2$ the partition cliques are the matching edges, and the cross cliques are the non-matching
edges, recovering the construction of \cref{defn:witness-lg}.

\begin{prop}[$G^*$ realizes $K_n$]\label{prop:witness-realizes}
    $G^*$ satisfies (i) and (ii) of \cref{lem:realization-rs}, so
    $\gex_{g^*}(G^*)\cong K_n$.
\end{prop}

\begin{proof}
    For (i), fix $v\in[n]$, and say $v\in P$. A clique
    of $\Omega$ through $v$ is either the partition clique $P$ or some cross clique $Q$ through $v$, each of which meets $P$ exactly at $v$. Thus
    $G^*[C_v]$ is a star with center $P$.

    For (ii), each pair $\{u,v\}\in\binom{[n]}{2}$ lies in some partition clique $P$ or else in a cross clique $Q$, and thus all such pairs are covered by $\Omega$.
\end{proof}

\begin{lem}[Treewidth of the witness graph]\label{lem:witness-tw}
    For $n=rm$ with $m\ge r+1$,
    \[
        \tw(G^*)=m-1
    \]
\end{lem}

\begin{proof}
    For the upper bound, as in \cref{lem:witness}, take a root bag
    $B_0=\mathcal P$ and a leaf bag $B_Q=\{Q,P_{i_1},\hdots,P_{i_r}\}$ where $Q\cap P_{i_j}\neq \emptyset$. Every clique of $\Omega$ lies in a bag; every edge $PQ$ lies in $B_Q$; and
    the bags containing a fixed vertex form a subtree. The sizes are $|B_0|=m$ and $|B_Q|=r+1$, so the width is
    $\max(m,r+1)-1=m-1$ since $m\ge r+1$.

    For the lower bound, \cref{lem:tw-growth} with $\mu=r$ gives
    $\tw(\gex_{g^*}(G^*))\le\max\{r\,\tw(G^*)+(r-1),\,2r-2\}$, since gadgets of $\LL_r$ have at most $2r-1$ vertices. As $n-1\ge r^2+r-1>2r-2$, $\gex_{g^*}(G^*)\cong
    K_n$ forces $\tw(G^*)\ge(n-r)/r=m-1$. The two bounds
    coincide.
\end{proof}

\begin{rem}\label{rem:small-corner}
    The width is determined by the root bag $B_0$ once $m\ge r+1$. For $m=r$ the leaf bags dominate and the decomposition gives
    only $\tw(G^*)\le r$; indeed, $P\cap Q\neq \emptyset$ implies $PQ\in E(G^*)$ for every $(P,Q)\in\mathcal P\times\mathcal Q$, so $G^*$ contains a complete bipartite $K_{r,|\mathcal Q|}$ of treewidth
    $r$, one above the lower bound $m-1$. In this case a different witness would be needed to achieve the lower bound where possible.
\end{rem}

\begin{lem}[Coefficient of the witness graph]\label{lem:witness-coeff}
    $\beta(G^*,K_n)\ge 1$; in particular $K_n$ occurs with nonzero coefficient in the
    gadget expansion of $G^*$.
\end{lem}

\begin{proof}
    Fix a realizing assignment $g\colon E(G^*)\to\LL_r$ and an isomorphism
    $\gex_g(G^*)\cong K_n$. Let $R$ be the identification graph on the vertices of the blow-up $B_r(G^*)$ whose
    edges connect vertices identified by the gadgets of $g$, i.e., the connected components of $R$ are the
    classes $[n]$ of $\sim_g$. Observe that $G^*$ admits the bipartition $\mathcal P,\mathcal Q$ into partition cliques and cross cliques. For $W\in\Omega$ write
    $\pi_g(W)\subseteq V(K_n)$ for the image of the vertices of $B_r(W)$; $|\pi_g(W)|=r$ since $B_r(W)$ is a
    clique and $K_n$ is loopless.

    First we claim that the images of the partition cliques $\{\pi_g(P)\}_{P\in\mathcal P}$ partition $V(K_n)$. Since $\sum_{P\in\mathcal P}|\pi_g(P)|=rm=n=|V(K_n)|$, it
    suffices to show they cover $V(K_n)$. Suppose for contradiction that some vertex $v\in K_n$ lies in no $\pi_g(P)$.
    Then the class $\pi_g^{-1}(v)\subseteq V(B_r(G^*))$ contains no vertex from a partition clique, and thus each edge incident to $v$ is an edge of some cross clique $Q\in\mathcal Q$. Since $v$ is adjacent to all other vertices of $K_n$ this implies

    \[
        V(K_n)\setminus\{v\}\subseteq\bigcup_{Q\;:\;v\,\in\,\pi_g(Q)}\bigl(\pi_g(Q)\setminus\{v\}\bigr).
    \]

    \noi Now each cross clique $Q$ contributes at most $r-1$ distinct vertices to the union, so there must be at least two cross cliques $Q_i\neq Q_j$ containing $v$. Thus $\pi_g^{-1}(v)$ contains a vertex $v_i\in V(B_r(Q_i))$ and a vertex $v_j\in V(B_r(Q_j))$. Now the bipartition $\mathcal P,\mathcal Q$ of $G^*$ implies each edge of the identification graph $R$ connects a partition clique to a cross clique. Thus there exists some path from $v_i$ to $v_j$ in $R$ which must contain a vertex of $B_r(P)$ for some partition clique $P\in\mathcal P$, the desired contradiction. Therefore
    $\{\pi_g(P)\}_{P\in\mathcal P}$ partitions $V(K_n)$.

    Now we claim that every realizing assignment uses only single gadgets, i.e., gadgets $\boldsymbol{F}_M\in\LL_r$ with $|M|=1$.
    Fix a cross clique $Q$ with neighborhood $\{P_{1},\hdots,P_r\}$ in $G^*$.
    Each gadget on $QP_i$ matches a nonempty subset of $B_r(Q)$ to $B_r(P_i)$.
    These $r$ subsets are pairwise disjoint, since the images $\pi_g(P_i)$ are disjoint.
    As $B_r(Q)$ has $r$ vertices, each subset is a singleton.
    Thus every edge incident to $Q$ is assigned a single gadget.

    Every realizing assignment therefore contributes
    $\prod_{e\in E(G^*)}\alpha_1=(+1)^{|E(G^*)|}=+1$, using $\alpha_1=+1$
    (\cref{lem:coefficients}). The canonical assignment $g^*$ is a realizer
    (\cref{prop:witness-realizes}), so
    \[
        \beta(G^*,K_n)
        =\sum_{g:\,\gex_g(G^*)\cong K_n}\ \prod_{e\in E(G^*)}\alpha(g(e))
        =\#\{\text{realizing assignments}\}\ge 1. \qedhere
    \]
\end{proof}

\begin{cor}[Transfer for the $(r,S)$-line graph]\label{cor:tightness}
    For every $k\ge1$,
    \[
        H_1\WL{rk+r-1}H_2\;\Longrightarrow\;\Lr(H_1)\WL{k}\Lr(H_2).
    \]
    Further, for $k\ge r$, $rk+r-1$ is the least possible value on the left-hand side: there exist $H_1\WL{rk+r-2}H_2$ with
    $\Lr(H_1)\not\WL{k}\Lr(H_2)$. The first statement holds for every $S\subseteq\{0,\dots,r-1\}$, the second whenever
    $\{1\}\subseteq S\subseteq\{1,\dots,r-1\}$.
\end{cor}

\begin{proof}
    For the upper bound, as in \cref{lem:upper-bound} with $\mu=r$: for $G\in\T_k$ and any
    $g\colon E(G)\to\LL_r$, \cref{lem:tw-growth} gives $\tw(\gex_g(G))\le rk+(r-1)$ (gadgets of $\LL_r$ have at most $2r$ vertices, and $2r-1\le rk+r-1$),
    so the isomorphism types $\Gamma_G$ of the loopless expansions $\gex_g(G)$ lie in $\T_{rk+(r-1)}$.\footnote
    {
        Without loss of generality $G$ has no isolated vertices, as
        \cref{defn:represents} requires: writing $G=G'\sqcup\overline K_i$ with $G'$ free of
        isolated vertices, $\hom(G,X)=|V(X)|^{\,i}\hom(G',X)$, and
        $|V(\Lr(H))|=\hom(K_r,H)/r!$ is a $\WLp{r-1}$-invariant, since every homomorphism from $K_r$
        into a loopless graph is injective.
    }
    By \eqref{eqn:decomp-iso} (with $d=r!$, \cref{lem:duplication-rs}) $\hom(G,\Lr(H))$ is a
    combination of the counts $\hom(\gamma,H)$, $\gamma\in\Gamma_G$; so
    $H_1\WL{rk+(r-1)}H_2$ forces $\Lr(H_1)\WL{k}\Lr(H_2)$.

    For the lower bound, set $n=r(k+1)$, $m=k+1\ge r+1$, and take the witness graph $G^*$ on
    $[n]$. By \cref{lem:witness-tw} $\tw(G^*)=k$; by \cref{prop:witness-realizes}
    $\gex_{g^*}(G^*)\cong K_n$; and by \cref{lem:witness-coeff}
    $\beta(G^*,K_n)\geq1$. Applying \eqref{eqn:decomp-iso} to $G^*$,
    \[
        \hom(G^*,\Lr(H))
        =(r!)^{-|V(G^*)|}\sum_{\gamma\in\Gamma}\beta(G^*,\gamma)\,\hom(\gamma,H),
        \qquad \Gamma=\{\gamma\in\G:\beta(G^*,\gamma)\neq0\}.
    \]
    Since $\beta(G^*,K_n)\geq1>0$, $K_n\in\Gamma$, and
    $\tw(K_n)=n-1=rk+(r-1)>rk+r-2$, so $\Gamma\not\subseteq\T_{rk+r-2}$. Were
    $\hom(G^*,\Lr(\cdot))$ a $\WL{rk+r-2}$-invariant, \cref{thm:hdc} applied to
    $(\Gamma,\beta(G^*,\cdot))$ would force $\Gamma\subseteq\T_{rk+r-2}$, a
    contradiction. Hence there exist $H_1\WL{rk+r-2}H_2$ with
    $\hom(G^*,\Lr(H_1))\neq\hom(G^*,\Lr(H_2))$; as $\tw(G^*)=k$, this gives
    $\Lr(H_1)\not\WL{k}\Lr(H_2)$. Together with the upper bound this is the
    case $k\ge r$ of \cref{thm:main-rs} where $\tau_{r,S}(k)=rk+r-1$.
\end{proof}

\begin{rem}\label{rem:fused}
    For $S=\{0,1,\dots,r-1\}$, $\LrS(H)$ is the complete graph on the $r$-cliques of $H$.
    Its order is $\hom(K_r,H)/r!$, a $\WLp{r-1}$-invariant since $\tw(K_r)=r-1$.
    Thus, for every $k\geq1$,
    \[
        H_1\WL{r-1}H_2\;\Longrightarrow\;\LrS(H_1)\WL{k}\LrS(H_2),
    \]
    so together with the lower bound of \cref{rem:k1-rs}, $\tau_{r,S}(k)=r-1$ for every $k\geq1$, and the lower bound in item~1 of \cref{thm:main-rs} is attained.
\end{rem}

\subsection{The case \texorpdfstring{$k=1$}{k = 1}}\label{subsec:k1-rs}

The witness graph of \cref{defn:witness-rs} has treewidth at least $r\ge2$ (\cref{lem:witness-tw,rem:small-corner}), and thus fails for $k=1$. Additionally, \cref{lem:tw-growth} gives the loose upper bound $\tau_{r,S}(1)\le 2r-1$, rather than the true value of $r-1$. The gap comes from the use of the blow-up $G\boxtimes K_\mu$, which fully connects adjacent cliques. Gadgets of $\LL_r$ add no edges between cliques, so the treewidth grows less.

\begin{lem}[Forest expansion]\label{lem:forest-rs}
    For every forest $G$ and assignment $g\colon E(G)\to\LL_r$, we have
    $\tw(\gex_g(G))\le r-1$.
\end{lem}
\begin{proof}
    We may assume that $G$ is a tree with at least one edge, since expansions respect disjoint unions and an isolated vertex expands to $r$ isolated vertices.
    Write $\gex_g(G)=B_r(G)/\!\sim_g$ as in \cref{defn:br}, with quotient map $\pi_g$, and
    let $R$ be the identification graph on $V(G)\times[r]$ whose edges record the
    identifications prescribed by the gadgets, so the components of $R$ are the vertices
    of $\gex_g(G)$. The proof of \cref{lem:R-forest} applies with $[r]$ in place of
    $[2]$: $R$ is a forest and the $r$ interface vertices of each $x\in V(G)$ lie in
    distinct components.

    Take $G$ itself as the decomposition tree, with bag
    $B_x=\pi_g(\{x\}\times[r])$ at $x\in V(G)$; by the above $|B_x|=r$. Every vertex of
    $\gex_g(G)$ is the class of some interface vertex and so lies in some bag. Every edge
    of $\gex_g(G)$ is the image of an edge of some $B_r(x)$, since gadgets carry no
    cross-interface edges, and so lies in $B_x$. Finally, for a fixed class $c$, the set
    $\{x: c\in B_x\}$ is the projection to $V(G)$ of a component of $R$, which is connected
    in $G$ because $R$ carries edges only along edges of $G$; a connected subgraph of a
    forest is a subtree. Thus $(G,(B_x)_x)$ is a tree decomposition of width $r-1$.
\end{proof}

\noindent At $r=2$ this recovers \cref{prop:forest-comp}: a graph of treewidth at most $1$
is a forest.

\begin{cor}[Transfer at $k=1$]\label{cor:k1-rs}
    $\tau_{r,S}(1)=r-1$ for every $r\ge2$ and $S\subseteq\{0,\dots,r-1\}$;
    this is the case $k=1$ of \cref{thm:main-rs}.
\end{cor}
\begin{proof}
    For the upper bound, write a forest $G$ as $G'\sqcup\overline K_i$, where $G'$ has no isolated vertices.
    Then $\hom(G,\Lr(H))=(\hom(K_r,H)/r!)^i\hom(G',\Lr(H))$.
    By \eqref{eqn:decomp-iso}, the last factor is
    a rational combination of the $\hom(\gamma,H)$ with $\gamma=\gex_g(G')\in\T_{r-1}$ by
    \cref{lem:forest-rs}.
    Since $K_r\in\T_{r-1}$, all these counts agree when $H_1\WL{r-1}H_2$, so
    $\Lr(H_1)\WL{1}\Lr(H_2)$.

    For the lower bound, for $r=2$ this is \cref{thm:main-lg}, so let $r\ge3$. Since every
    homomorphism from $K_r$ into a loopless graph is injective, with image an $r$-clique,
    $|V(\Lr(H))|=\hom(K_r,H)/r!$. Were this a $\WL{r-2}$-invariant, \cref{thm:hdc} would
    place $K_r$ in $\T_{r-2}$, contradicting $\tw(K_r)=r-1$. Hence some $H_1\WL{r-2}H_2$ have
    $|V(\Lr(H_1))|\ne|V(\Lr(H_2))|$, i.e., $\hom(K_1,\Lr(H_1))\ne\hom(K_1,\Lr(H_2))$, and as
    $K_1\in\T_1$ this gives $\Lr(H_1)\not\WL{1}\Lr(H_2)$.
\end{proof}

\begin{rem}\label{rem:k1-rs}
    The lower bound uses only the probe $K_1$, so $\tau_{r,S}(k)\ge r-1$ for every $k\ge1$;
    this is the bound that \cref{cor:cycle-lower} doubles for $k\ge2$. The upper bound uses
    only that gadgets add no cross-interface edges. For $r\ge3$ there is no analogue of the
    covering-theoretic argument available at $r=2$: a cover $X\to H$ need not induce a
    cover $\Lr(X)\to\Lr(H)$, since a triangle of $H$ lifts to a path.
\end{rem}

\subsection{Lower bound for \texorpdfstring{$2\leq k\leq r-1$}{2 <= k <= r-1}}\label{subsec:cycle}
\cref{cor:tightness} leaves open sharpness in the range $2\le k\le r-1$, where the witness graphs of \cref{defn:witness-rs}, all of treewidth at least $r>k$ (\cref{lem:witness-tw,rem:small-corner}), are too dense. We record a lower bound that holds uniformly in $k\ge2$, obtained from the alternative witness graph $C_{2r-1}$ of treewidth $2$. Throughout this
subsection we assume $\{1\}\subseteq S\subseteq\{1,\dots,r-1\}$ and the sign condition
\begin{equation}\label{eqn:dagger}\tag{$\dagger$}
    (-1)^{t-1}\alpha_t\ \ge\ 0\quad(1\le t\le r),\qquad \alpha_{r-1}\neq0.
\end{equation}
The following gives some criteria under which this holds automatically.

\begin{lem}[Admissible intersection sets]\label{lem:dagger-suff}
    Write $h_S(t)=(-1)^{t-1}\alpha_t=\sum_{s\in S}(-1)^{s-1}\binom ts$, so that
    \eqref{eqn:dagger} reads $h_S(t)\ge0$ for $1\le t\le r$. Then \eqref{eqn:dagger}
    holds in both of the following cases:
    \begin{enumerate}
        \item[(i)] $S$ consists of only odd numbers, in which case
              $h_S(t)=\sum_{s\in S}\binom ts\ge0$ for every $t$;
        \item[(ii)] $S=\{1,\dots,r-1\}$, in which case $h_S(t)=1+(-1)^{r}\binom{t-1}{r-1}$,
              which is nonnegative for $1\le t\le r$.
    \end{enumerate}
    In particular \eqref{eqn:dagger} holds for $S=\{1\}$
    and for the hypergraph-line-graph family $S=\{1,\dots,r-1\}$.
\end{lem}

\begin{proof}
    We treat the two conditions separately.
    For (i), when every $s\in S$ is odd then $h_S(t)$ is a sum of binomial coefficients.

    For (ii), for $S=\{1,\dots,r-1\}$, split off the $s=0$ term and apply $\sum_{s=0}^{r-1}(-1)^s\binom ts=(-1)^{r-1}\binom{t-1}{r-1}$:
        \[
            h_S(t)=\sum_{s=1}^{r-1}(-1)^{s-1}\binom ts
            =1-\sum_{s=0}^{r-1}(-1)^{s}\binom ts
            =1-(-1)^{r-1}\binom{t-1}{r-1}
            =1+(-1)^{r}\binom{t-1}{r-1}.
        \]
        For $t<r$ the binomial vanishes and $h_S(t)=1$; for $t=r$ it is
        $1+(-1)^{r}\in\{0,2\}$.
        In both cases $\alpha_{r-1}\neq0$, since $h_S(r-1)\neq0$: in case (i) the term
        $s=1$ alone gives $h_S(r-1)\ge\binom{r-1}{1}>0$, and in case (ii) $h_S(r-1)=1$.
\end{proof}

For $p\ge3$, as in \cref{lem:witness-coeff}, a realizing assignment $g\colon E(C_p)\to\LL_r$ with
$\gex_g(C_p)\cong K_n$ partitions the $rp$ vertices of $B_r(C_p)$ into the $n$ classes
$\pi_g^{-1}(1),\dots,\pi_g^{-1}(n)$, which are the components of the identification graph
$R$. We define the \emph{support} $V_i\subseteq V(C_p)$ of the class $\pi_g^{-1}(i)$ as
the set of vertices $v\in V(C_p)$ whose blow-up $B_r(v)$ meets $\pi_g^{-1}(i)$.

\begin{lem}[Supports]\label{lem:supports}
    Let $g$ be a realizing assignment with $\gex_g(C_p)\cong K_n$ and $V_1,\dots,V_n$ the
    supports of its classes.
    \begin{enumerate}
        \item[(i)] Each vertex $v$ of $C_p$ lies in exactly $r$ supports, and $\pi_g^{-1}(i)\to V_i$ is a bijection for all $i$.
        \item[(ii)] Every edge of $R$ joins $B_r(u)$ to $B_r(v)$ for some edge $uv$ of
            $C_p$. Thus each $V_i$ is connected in $C_p$, i.e., an \emph{arc}
            (vertex set inducing a path) or all of $V(C_p)$.
        \item[(iii)] $\gex_g(C_p)$ is complete iff $V_1,\dots,V_n$ intersect pairwise.
    \end{enumerate}
\end{lem}
\begin{proof}
    For (i), the restriction $\pi_g\colon B_r(v)\to K_n$ is injective, since $B_r(v)$ is a
    clique and $K_n$ is loopless. Thus each class $\pi_g^{-1}(i)$ intersects $B_r(v)$ in at
    most one vertex, and the $r$ vertices of $B_r(v)$ lie in $r$ distinct
    classes. This implies each vertex $v\in V(C_p)$ lies in exactly $r$ supports, and the map $\pi_g^{-1}(i)\to V(C_p)$ sending vertices of $B_r(v)$ to $v$ is injective; by the definition of the support, its image is $V_i$.

    For (ii), the edges of $R$ are the identifications prescribed by the gadgets, and $g(uv)$
    identifies vertices of $B_r(u)$ with vertices of $B_r(v)$ only. Mapping the
    component $\pi_g^{-1}(i)$ of $R$ to $C_p$ therefore gives a connected subgraph of
    $C_p$ with vertex set $V_i$.

    For (iii), two classes are adjacent in $\gex_g(C_p)$ iff some $B_r(v)$ contains a vertex of
    each, i.e., iff $v\in V_i\cap V_j$ for some $v$.
\end{proof}

\begin{lem}[Arc packing]\label{lem:arc-packing}
    Let $\mathcal A$ be a family of $m$ arcs of $C_p$ that intersect pairwise and
    cover each vertex exactly $d$ times. Then $m\le 2d-1$.
\end{lem}
\begin{proof}
    Let $A\in\mathcal A$ with endpoints $a,b$ be an arc of minimal length. Any other arc
    of $\mathcal A$ intersects $A$ and must contain one of its endpoints by the minimality of $A$. Since each vertex is covered exactly $d$ times, at most $d$ arcs contain $a$ and at most $d$ arcs contain $b$. Since $A$ contains both, the total number of arcs $m\leq 2d-1$.
\end{proof}

\begin{prop}[Cycle witness]\label{prop:cycle-witness}
    Let $r\ge3$ and $n=2r-1$. Then for every $S$ satisfying \eqref{eqn:dagger},
    \[
        (-1)^r\,\beta\bigl(C_{n},\,K_{n}\bigr)\;\ge\;1,
    \]
    where
    \[
        \beta(C_{n},K_{n})=\sum_{g\,:\,\gex_g(C_{n})\cong K_{n}}\;
            \prod_{e\in E(C_{n})}\alpha(g(e))
    \]
    is the coefficient in the gadget decomposition of \cref{lem:gadget-decomp}.
\end{prop}
\begin{proof}
    Fix a realizing assignment $g$ with $\gex_g(C_n)\cong K_n$, and let
    $V_1,\dots,V_n$ be the supports of its classes.

    We first show that no support is all of $V(C_n)$.
    Suppose $j\ge1$ supports equal $V(C_n)$. By
    \cref{lem:supports}, the remaining $n-j$ supports are arcs, intersecting pairwise and
    covering each vertex of $C_n$ exactly $r-j$ times. \cref{lem:arc-packing} gives
    $n-j\le2(r-j)-1$, i.e., $2r-1-j\le2r-2j-1$, i.e., $j\le0$, a contradiction. So every
    $V_i$ is an arc, and $C_n[V_i]$ is a path.

    Next, by results (i) and (ii) of \cref{lem:supports}, the component
    of $R$ on $\pi_g^{-1}(i)$ is isomorphic to
    $C_n[V_i]$, i.e., it is a path on $|V_i|$ vertices. For an edge $e$ of $C_n$ let
    $t_e$ be the number of supports containing $e$. A class whose support contains
    $e=uv$ has one vertex in each of $B_r(u)$ and $B_r(v)$, joined by exactly one edge
    of $R$. So $g(e)=\boldsymbol F_M$ for a matching $M$ with $|M|=t_e$, and by \cref{defn:family-rs}
    $\alpha(g(e))=\alpha_{t_e}$, so the summand of $g$ in
    $\beta(C_n,K_n)$ is $\prod_e\alpha_{t_e}$.
    Under \eqref{eqn:dagger} every nonzero $\alpha_{t_e}$ has sign $(-1)^{t_e-1}$, so
    every nonzero summand of $\beta(C_n,K_n)$ has sign $(-1)^{\sum_e(t_e-1)}$. Using \cref{lem:supports}(i) to count the edges of $R$ in two
    ways, we have
    \[
        \sum_{e\in E(C_n)}t_e=\sum_{i}(|V_i|-1)=rn-n=(r-1)n.
    \]
    Thus the sign of the summand of $g$ is
    \[
        \prod_{e\in E(C_n)}(-1)^{t_e-1}=(-1)^{\sum_e(t_e-1)}=(-1)^{(r-1)n-n}=(-1)^{(r-2)n}=(-1)^r
    \]
    the last equality because $n$ is odd, and therefore the summands cannot cancel.

    Finally, we exhibit a realizing assignment. Identify $V(C_n)$ with $\Z/n\Z$, orient each edge $e_v=\{v,v+1\}$ from $v$ to $v+1$, and label the $r$ interface vertices of $B_r(v)$ bijectively by $\{v-r+1,\dots,v\}\subseteq\Z/n\Z$, a set of $r$ distinct residues since $r<n$.
    Let $g(e_v)=\boldsymbol F_M$, where $M$ identifies the unique vertices of $B_r(v)$ and $B_r(v+1)$ of label $i$ for each of the $r-1$ labels $i\in\{v-r+2,\dots,v\}$ common to both blocks. Then $|M|=r-1$, and $\alpha_{r-1}\neq0$ by \eqref{eqn:dagger}, so $g(e_v)\in\LL_r$.
    Observe that label $i$ appears exactly in the blocks $B_r(v)$ for $v$ in the arc $V_i=\{i,\dots,i+r-1\}$, and $g$ identifies only vertices of equal label, gluing consecutive blocks along each arc; hence the classes of $\sim_g$ are exactly the $n$ label classes, with supports $V_i$. The vertices of each $B_r(v)$ have distinct labels, so $\gex_g(C_n)$ is loopless. Since any two arcs of $r$ vertices meet when $2r>n$, \cref{lem:supports}(iii) shows that $\gex_g(C_n)$ is complete.
    Finally, each edge $e_v$ lies in exactly $r-1$ supports, namely $V_i$ for $i\in\{v-r+2,\dots,v\}$, so $t_{e_v}=r-1$ throughout and the summand of $g$ is $\alpha_{r-1}^{\,n}\neq0$. Hence
    $(-1)^r\beta(C_n,K_n)\ge1$.
\end{proof}

\begin{cor}[Uniform lower bound for $2\le k\le r-1$]\label{cor:cycle-lower}
    Let $r\ge3$ and let $S$ satisfy \eqref{eqn:dagger}. For every $k\ge2$ there exist
    $H_1\WL{2r-3}H_2$ with $\Lr(H_1)\not\WL{k}\Lr(H_2)$; equivalently,
    $\tau_{r,S}(k)\ge 2r-2$.
\end{cor}
\begin{proof}
    $C_{2r-1}$ has treewidth $2\le k$. By \cref{prop:cycle-witness}
    $\beta(C_{2r-1},K_{2r-1})\ne0$, so $K_{2r-1}\in\Gamma=\{\gamma\in\G:\beta(C_{2r-1},\gamma)\ne0\}$, and
    $\tw(K_{2r-1})=2r-2$. Were $\hom(C_{2r-1},\Lr(\cdot))$ a $\WL{2r-3}$-invariant, then
    \cref{thm:hdc} applied to $(\Gamma,\beta(C_{2r-1},\cdot))$ through \eqref{eqn:decomp-iso}
    would force $\Gamma\subseteq\T_{2r-3}$, contradicting $2r-2>2r-3$. Hence some
    $H_1\WL{2r-3}H_2$ satisfy $\hom(C_{2r-1},\Lr(H_1))\ne\hom(C_{2r-1},\Lr(H_2))$; as
    $\tw(C_{2r-1})=2\le k$, this gives $\Lr(H_1)\not\WL{k}\Lr(H_2)$.
\end{proof}

\begin{rem}\label{rem:cycle-scope}
    This doubles the bound $\tau_{r,S}(k)\ge r-1$ of \cref{rem:k1-rs}; for $r=3$ it gives
    $\tau_{3,S}(2)\ge4$. By \cref{lem:supports} and \cref{lem:arc-packing} no cycle realizes $K_n$ with $n>2r-1$, so $2r-2$ is the best
    bound available in this way.
    For $\{1\}\subseteq S\subseteq\{1,\dots,r-1\}$, it remains to determine whether $\tau_{r,S}(k)\ge2r-2$ holds for $2\le k\le r-1$ without the sign condition \eqref{eqn:dagger}.
\end{rem}

\section*{Use of AI}
We used Claude Fable 5.1 and Opus 5.5 throughout this project. The idea of using gadget expansions to study line graphs and WL-equivalence was our own. Claude assisted with literature review and wrote code to test intermediate claims computationally. The witness graph constructions of our lower bounds are based on small examples found in searches written and executed by Claude. Claude produced the figures, checked the proofs, and proofread the manuscript. Claude also drafted parts of the text, most of which were completely rewritten by the authors. The authors have verified the paper's contents and are responsible for any errors.

\bibliographystyle{alpha}
\bibliography{references}

\newcommand{\etalchar}[1]{$^{#1}$}
\begin{thebibliography}{MRF{\etalchar{+}}19}

\bibitem[ADW17]{ADW17}
Samson Abramsky, Anuj Dawar, and Pengming Wang.
\newblock The pebbling comonad in finite model theory.
\newblock In {\em 32nd Annual ACM/IEEE Symposium on Logic in Computer Science (LICS)}, pages 1--12, 2017.

\bibitem[AMP16]{AMP16}
Kwangjun Ahn, Dhruv Medarametla, and Aaron Potechin.
\newblock Graph matrices: norm bounds and applications.
\newblock {\em arXiv preprint arXiv:1604.03423}, 2016.

\bibitem[Bab16]{Babai16}
L{\'a}szl{\'o} Babai.
\newblock Graph isomorphism in quasipolynomial time [extended abstract].
\newblock In {\em Proceedings of the 48th Annual ACM SIGACT Symposium on Theory of Computing (STOC)}, pages 684--697, 2016.

\bibitem[BB17]{BB17}
Jacob Biamonte and Ville Bergholm.
\newblock Tensor networks in a nutshell.
\newblock {\em arXiv preprint arXiv:1708.00006}, 2017.

\bibitem[BCN89]{BCN89}
A.~E. Brouwer, A.~M. Cohen, and A.~Neumaier.
\newblock {\em Distance-Regular Graphs}.
\newblock Springer, Berlin, 1989.

\bibitem[Ber73]{Berge73}
Claude Berge.
\newblock {\em Graphs and Hypergraphs}.
\newblock North-Holland, Amsterdam, 1973.

\bibitem[BES80]{BES80}
L{\'a}szl{\'o} Babai, Paul Erd{\H o}s, and Stanley~M. Selkow.
\newblock Random graph isomorphism.
\newblock {\em SIAM Journal on Computing}, 9(3):628--635, 1980.

\bibitem[BHK{\etalchar{+}}19]{BHK19}
Boaz Barak, Samuel~B. Hopkins, Jonathan Kelner, Pravesh~K. Kothari, Ankur Moitra, and Aaron Potechin.
\newblock A nearly tight sum-of-squares lower bound for the planted clique problem.
\newblock {\em SIAM Journal on Computing}, 48(2):687--735, 2019.

\bibitem[BI84]{BannaiIto84}
Eiichi Bannai and Tatsuro Ito.
\newblock {\em Algebraic Combinatorics {I}: Association Schemes}.
\newblock Benjamin/Cummings, Menlo Park, CA, 1984.

\bibitem[CDM17]{CDM17}
Radu Curticapean, Holger Dell, and D{\'a}niel Marx.
\newblock Homomorphisms are a good basis for counting small subgraphs.
\newblock In {\em Proceedings of the 49th Annual ACM SIGACT Symposium on Theory of Computing}, pages 210--223, 2017.

\bibitem[CFI92]{CFI92}
Jin-Yi Cai, Martin F{\"u}rer, and Neil Immerman.
\newblock An optimal lower bound on the number of variables for graph identification.
\newblock {\em Combinatorica}, 12(4):389--410, 1992.

\bibitem[CL87]{ChenLih87}
Bor-Liang Chen and Ko-Wei Lih.
\newblock Hamiltonian uniform subset graphs.
\newblock {\em Journal of Combinatorial Theory, Series B}, 42(3):257--263, 1987.

\bibitem[CLB19]{ChenLiBruna19}
Zhengdao Chen, Xiang Li, and Joan Bruna.
\newblock Supervised community detection with line graph neural networks.
\newblock In {\em International Conference on Learning Representations (ICLR)}, 2019.

\bibitem[DGR18]{DGR18}
Holger Dell, Martin Grohe, and Gaurav Rattan.
\newblock Lov{\'a}sz meets {Weisfeiler} and {Leman}.
\newblock In {\em 45th International Colloquium on Automata, Languages, and Programming (ICALP 2018)}, volume 107 of {\em LIPIcs}, pages 40:1--40:14. Schloss Dagstuhl, 2018.

\bibitem[DHP24]{DHP24}
Amina Doumane, Samuel Humeau, and Damien Pous.
\newblock A finite presentation of graphs of treewidth at most three.
\newblock In {\em 51st International Colloquium on Automata, Languages, and Programming (ICALP 2024)}, volume 297 of {\em LIPIcs}, pages 135:1--135:18. Schloss Dagstuhl, 2024.

\bibitem[DJR21]{DJR21}
Anuj Dawar, Tom{\'a}{\v s} Jakl, and Luca Reggio.
\newblock Lov{\'a}sz-type theorems and game comonads.
\newblock In {\em 36th Annual ACM/IEEE Symposium on Logic in Computer Science (LICS)}, pages 1--13, 2021.

\bibitem[Dvo10]{Dvorak10}
Zden{\v{e}}k Dvo{\v{r}}{\'a}k.
\newblock On recognizing graphs by numbers of homomorphisms.
\newblock {\em Journal of Graph Theory}, 64(4):330--342, 2010.

\bibitem[FSS24]{FSS24}
Eva Fluck, Tim Seppelt, and Gian~Luca Spitzer.
\newblock Going deep and going wide: Counting logic and homomorphism indistinguishability over graphs of bounded treedepth and treewidth.
\newblock In {\em 32nd EACSL Annual Conference on Computer Science Logic (CSL 2024)}, volume 288 of {\em LIPIcs}, pages 27:1--27:17. Schloss Dagstuhl, 2024.

\bibitem[GN21]{GroheNeuen21}
Martin Grohe and Daniel Neuen.
\newblock Recent advances on the graph isomorphism problem.
\newblock In {\em Surveys in Combinatorics 2021}, London Mathematical Society Lecture Note Series, pages 187--234. Cambridge University Press, 2021.

\bibitem[GRS25]{GRS22}
Martin Grohe, Gaurav Rattan, and Tim Seppelt.
\newblock Homomorphism tensors and linear equations.
\newblock {\em Advances in Combinatorics}, April 2025.

\bibitem[Hav18]{Haviv18}
Ishay Haviv.
\newblock On minrank and the {L}ov\'asz theta function.
\newblock In {\em Approximation, Randomization, and Combinatorial Optimization (APPROX/RANDOM 2018)}, volume 116 of {\em LIPIcs}, pages 13:1--13:15, 2018.

\bibitem[Kra43]{Krausz43}
J.~Krausz.
\newblock D{\'e}monstration nouvelle d'un th{\'e}or{\`e}me de {Whitney} sur les r{\'e}seaux.
\newblock {\em Matematikai {\'e}s Fizikai Lapok}, 50:75--85, 1943.

\bibitem[KV15]{KV15}
Andreas Krebs and Oleg Verbitsky.
\newblock Universal covers, color refinement, and two-variable counting logic: Lower bounds for the depth.
\newblock In {\em Proceedings of the 30th Annual ACM/IEEE Symposium on Logic in Computer Science (LICS 2015)}, pages 689--700. IEEE Computer Society, 2015.

\bibitem[Lei82]{Leighton82}
Frank~Thomson Leighton.
\newblock Finite common coverings of graphs.
\newblock {\em Journal of Combinatorial Theory, Series B}, 33(3):231--238, 1982.

\bibitem[Lov67]{Lovasz67}
L{\'a}szl{\'o} Lov{\'a}sz.
\newblock Operations with structures.
\newblock {\em Acta Mathematica Academiae Scientiarum Hungaricae}, 18:321--328, 1967.

\bibitem[Lov12]{Lovasz12}
L{\'a}szl{\'o} Lov{\'a}sz.
\newblock {\em Large Networks and Graph Limits}, volume~60 of {\em Colloquium Publications}.
\newblock American Mathematical Society, 2012.

\bibitem[LS22]{LiuSong22}
Xin Liu and Yangqiu Song.
\newblock Graph convolutional networks with dual message passing for subgraph isomorphism counting and matching.
\newblock {\em Proceedings of the AAAI Conference on Artificial Intelligence}, 36(7):7594--7602, 2022.

\bibitem[LSBT15]{GeneralizedLineGraphs15}
Aparna Lakshmanan~S., Csilla Bujt\'as, and Zsolt Tuza.
\newblock Generalized line graphs: {C}artesian products and complexity of recognition.
\newblock {\em The Electronic Journal of Combinatorics}, 22(3):P3.33, 2015.

\bibitem[LSBT16]{LakshmananBujtasTuza}
Aparna Lakshmanan~S., Csilla Bujt\'as, and Zsolt Tuza.
\newblock Induced cycles in triangle graphs.
\newblock {\em Discrete Applied Mathematics}, 209:264--275, 2016.

\bibitem[MP14]{McKayPiperno14}
Brendan~D. McKay and Adolfo Piperno.
\newblock Practical graph isomorphism, {II}.
\newblock {\em Journal of Symbolic Computation}, 60:94--112, 2014.

\bibitem[MR20]{MR20}
Laura Man{\v{c}}inska and David~E. Roberson.
\newblock Quantum isomorphism is equivalent to equality of homomorphism counts from planar graphs.
\newblock In {\em Proceedings of the 61st Annual IEEE Symposium on Foundations of Computer Science (FOCS)}, pages 661--672, 2020.
\newblock Full version: arXiv:1910.06958.

\bibitem[MRF{\etalchar{+}}19]{Morris19}
Christopher Morris, Martin Ritzert, Matthias Fey, William~L. Hamilton, Jan~Eric Lenssen, Gaurav Rattan, and Martin Grohe.
\newblock {Weisfeiler} and {Leman} go neural: Higher-order graph neural networks.
\newblock {\em Proceedings of the AAAI Conference on Artificial Intelligence}, 33(1):4602--4609, 2019.

\bibitem[MW19]{MW19}
Ankur Moitra and Alexander~S. Wein.
\newblock Spectral methods from tensor networks.
\newblock In {\em Proceedings of the 51st Annual ACM SIGACT Symposium on Theory of Computing (STOC 2019)}, pages 926--937, 2019.

\bibitem[Neu24]{Neuen23}
Daniel Neuen.
\newblock Homomorphism-distinguishing closedness for graphs of bounded tree-width.
\newblock In {\em 41st International Symposium on Theoretical Aspects of Computer Science (STACS 2024)}, volume 289 of {\em LIPIcs}, pages 53:1--53:12. Schloss Dagstuhl, 2024.

\bibitem[NRSS82]{NRSS82}
R.~N. Naik, S.~B. Rao, S.~S. Shrikhande, and N.~M. Singhi.
\newblock Intersection graphs of $k$-uniform linear hypergraphs.
\newblock {\em European Journal of Combinatorics}, 3(2):159--172, 1982.

\bibitem[NS26]{NS26}
Daniel Neuen and Tim Seppelt.
\newblock Distinguishing graphs by counting homomorphisms from sparse graphs.
\newblock In {\em 41st Annual Symposium on Logic in Computer Science (LICS 2026)}, volume 380 of {\em LIPIcs}, pages 75:1--75:29. Schloss Dagstuhl, 2026.

\bibitem[Ott97]{Otto97}
Martin Otto.
\newblock {\em Bounded Variable Logics and Counting: A Study in Finite Models}, volume~9 of {\em Lecture Notes in Logic}.
\newblock Springer, 1997.

\bibitem[Rob22]{Roberson22}
David~E. Roberson.
\newblock Oddomorphisms and homomorphism indistinguishability over graphs of bounded degree.
\newblock {\em arXiv preprint arXiv:2206.10321}, 2022.

\bibitem[RS23]{RobersonSeppelt23}
David~E. Roberson and Tim Seppelt.
\newblock {Lasserre} hierarchy for graph isomorphism and homomorphism indistinguishability.
\newblock In {\em 50th International Colloquium on Automata, Languages, and Programming (ICALP 2023)}, volume 261 of {\em LIPIcs}, pages 101:1--101:18, 2023.

\bibitem[Sep24a]{HomIndZoo}
Tim Seppelt.
\newblock Homomorphism indistinguishability zoo.
\newblock \url{https://tseppelt.github.io/homind-database}, 2024.
\newblock Online database.

\bibitem[Sep24b]{Seppelt23}
Tim Seppelt.
\newblock Logical equivalences, homomorphism indistinguishability, and forbidden minors.
\newblock {\em Information and Computation}, 301:105224, 2024.

\bibitem[Thi00]{Thi00}
Nicolas~M. Thi{\'e}ry.
\newblock Algebraic invariants of graphs; a study based on computer exploration.
\newblock {\em ACM SIGSAM Bulletin}, 34(3):9--20, 2000.

\bibitem[Wes01]{West2001}
Douglas~B. West.
\newblock {\em Introduction to Graph Theory}.
\newblock Prentice Hall, 2nd edition, 2001.

\bibitem[Whi32]{Whitney32}
Hassler Whitney.
\newblock Congruent graphs and the connectivity of graphs.
\newblock {\em American Journal of Mathematics}, 54(1):150--168, 1932.

\bibitem[XHLJ19]{Xu19}
Keyulu Xu, Weihua Hu, Jure Leskovec, and Stefanie Jegelka.
\newblock How powerful are graph neural networks?
\newblock In {\em International Conference on Learning Representations (ICLR)}, 2019.

\bibitem[Yan26]{Yang26}
Fan Yang.
\newblock On the expressive power of implicit line-graph higher-order {Weisfeiler}--{Leman}.
\newblock {\em arXiv preprint arXiv:2609.16412}, 2026.

\bibitem[YH24]{YangHuang24}
Fan Yang and Xingyue Huang.
\newblock Theoretical insights into line graph transformation on graph learning.
\newblock In {\em NeurIPS 2024 Workshop on Symmetry and Geometry in Neural Representations (NeurReps)}, 2024.

\end{thebibliography}

\end{document}